\documentclass[notitlepage,12pt]{article}
\usepackage{amsfonts,latexsym}

\title{{\Large Finite Difference Approximation \\
 of Free Discontinuity Problems}}
\author{
{\bf {\normalsize Massimo Gobbino}}\\ {\footnotesize Universit\`a degli 
Studi di Pisa}\vspace{-.1cm}
\\{\footnotesize Dipartimento di Matematica Applicata 
``Ulisse Dini''}\vspace{-.1cm}\\
{\footnotesize via Bonanno 25B, 56126
 PISA (Italy)}\vspace{-.1cm}\\
 {\footnotesize e-mail\vs\vs: \texttt{m.gobbino@dma.unipi.it}}
\and
{\bf {\normalsize Maria Giovanna Mora}}\\ {\footnotesize 
S.I.S.S.A.}\vspace{-.1cm}\\
{\footnotesize via Beirut 2/4, 34014 TRIESTE (Italy)}\vspace{-.1cm}\\
{\footnotesize e-mail: \texttt{mora@sissa.it}}}
\date{}

\newcommand{\ep}{\varepsilon}
\newcommand{\re}{{\mathbb{R}}}
\newcommand{\n}{{\mathbb{N}}}
\newcommand{\z}{{\mathbb{Z}}}

\newcommand{\lra}{\longrightarrow}
\newcommand{\nlra}{\rightarrow}

\newcommand{\vs}{\vspace{.2cm}}

\newcommand{\tc}{:\;}

\newcommand{\prf}{{\sc Proof.}$\;$}
\newcommand{\qed}{\hspace{1em}{$\Box$}\bigskip}
\newcommand{\qedbis}{\hspace{1em}{\Box}}
\newcommand{\sse}{\subseteq}

\newcommand{\foralll}{\forall\:}

\newcommand{\ie}{\emph{i.e.}}
\newcommand{\eg}{\emph{e.g.}}

\newcommand{\gconv}{$\Gamma$-convergence}

\newcommand{\hn}{{\cal H}^{n-1}}

\newcommand{\fep}{F_{\ep}}
\newcommand{\uep}{u_{\ep}}

\newcommand{\phiep}{\varphi_{\ep}}
\newcommand{\gmlim}[1]{\Gamma^{\mbox{-}}\!\mbox{-}\lim_{#1}}
\newcommand{\gmliminf}[1]{\Gamma^{\mbox{-}}\!\mbox{-}\liminf_{#1}}
\newcommand{\gmlimsup}[1]{\Gamma^{\mbox{-}}\!\mbox{-}\limsup_{#1}}

\newtheorem{thm}{Theorem}[section]
\newtheorem{rmk}[thm]{Remark}

\newtheorem{defn}[thm]{Definition}
\newtheorem{cor}[thm]{Corollary}
\newtheorem{example}[thm]{Example}
\newtheorem{lemma}[thm]{Lemma}

\begin{document}
\maketitle

\vspace{.6cm}
\begin{abstract}
\noindent{\footnotesize We approximate functionals depending on the gradient of $u$ and on the
behaviour of $u$ near the discontinuity points, by families of non-local
functionals where the gradient is replaced by finite differences. We
prove pointwise convergence, \gconv, and a compactness result which
implies, in particular, the convergence of minima and minimizers.}\vspace{.5cm}\\

\noindent{\footnotesize {\bf AMS (MOS) subject classifications:} 49J45,
65K10.}

\vspace{.5cm}
\noindent{\footnotesize {\bf Key words:} free discontinuity problems, $\mathit{SBV}$
functions, \gconv, non-local functionals.}
\thispagestyle{empty}
%

\end{abstract}

\begin{center}
\vspace{1.5cm} {\normalsize Ref. S.I.S.S.A. 14/99/M  (February 99)}
\end{center}
\pagebreak
\setcounter{page}{1}
\sloppy

\section{Introduction}
\label{sec:introduction}

In mathematical literature many free discontinuity problems have
been considered. The canonical examples are the minimum problems related
to the so called Mumford-Shah functional, defined by
\begin{equation}
	MS(u)=\int_{\Omega}|\nabla u(x)|^2\,dx+\hn(S_u),
	\label{defn:ms}
\end{equation}
where $\Omega$ is an open subset of $\re^n$, $u$ belongs to the space
$\mathit{SBV}(\Omega)$ of special functions with bounded variation (see
\S\,\ref{subsec:sbv}), $\nabla u$ is the approximate gradient of $u$,
$S_{u}$ is the set of essential discontinuity points of $u$, and $\hn$
is the $(n-1)$-dimensional Hausdorff measure.

This functional is the weak formulation in the space $\mathit{SBV}(\Omega)$ of
the functional introduced by {\sc D.  Mumford} and {\sc J.  Shah} in
\cite{ms} to approach image segmentation problems.

A natural generalization of (\ref{defn:ms}) are the functionals
\begin{equation}
	\mathcal{F}(u)=\int_{\Omega}\varphi(|\nabla u(x)|)\,dx+
	\int_{S_{u}}\psi(|u^{+}(x)-u^{-}(x)|)\,d\hn(x),
	\label{defn:fdp-gen}
\end{equation}
where $\varphi,\psi:[0,+\infty[\nlra [0,+\infty]$ are given
functions, and $u^{+}(x)$ and $u^{-}(x)$ are the approximate (in the measure
theoretic sense) $\limsup$ and $\liminf$ of $u$ at the point $x$.

By the semicontinuity and compactness theorem in $\mathit{SBV}$ proved by {\sc
L. Ambrosio} in \cite{ambrsc}, variational problems involving
$\mathcal{F}$ can be solved using the direct methods of the calculus
of variations: the interested reader can find appropriate references
in the survey \cite{ambrosiosurvey}.

Approximations of (\ref{defn:ms}) and (\ref{defn:fdp-gen}) have been
deeply studied in last years, both because of numerical applications,
and in order to approach evolution problems with free
discontinuities (cf. \cite{gobbino:msgf}). In this context,
approximation is always required in the sense of \gconv\ (see
\S\ \ref{subsec:gconv}), since this notion is stable under continuous
perturbations, and guarantees the convergence of minima and minimizers.

It is well known (cf.  \cite{bdm}) that functionals like
(\ref{defn:ms}) and (\ref{defn:fdp-gen}) \emph{cannot} be approximated
in the sense of \gconv\ by local integral functionals like
\begin{equation}
	\int_\Omega f_\ep(\nabla u(x))\,dx,
	\label{defn:loc-int-funct}
\end{equation}
defined in the Sobolev space $W^{1,2}(\Omega)$.
This difficulty has been overcome in different ways (cf. the survey
\cite{braides:survey}):
\begin{itemize}
	\item  by introducing an auxiliary variable as in \cite{at1,at2};

	\item by considering non-local functionals depending on the
	average of the gradient in small balls as in \cite{bdm};

	\item  by adding to (\ref{defn:loc-int-funct}) a singular
	perturbation depending on higher order derivatives of $u$
	(see \cite{stella1,stella2});

	\item  by using finite elements approximations, \ie\ local
	functionals like (\ref{defn:loc-int-funct}) defined in suitable
	spaces of piecewise affine functions (see \cite{chambolle:siam,cdm});

	\item  by considering non-local functionals where the gradient is
	replaced by finite differences (see \cite{gobbino:ms} and
	\cite{chambolle:num} for a numerical implementation).
\end{itemize}

The last approach was suggested in 1996 by \textsc{E. De Giorgi}, who
conjectured the convergence of the family
$$\mathcal{DG}_\ep(u)=\frac{1}{\ep}\int_{\re^n\times\re^n}\arctan\left(
\frac{(u(x+\ep\xi)-u(x))^2}{\ep}\right)e^{-|\xi|^2}\,d\xi\,dx, $$
to the Mumford-Shah functional in $\re^n$ (up to some constants),
both in the sense of pointwise convergence, and in the sense of
\gconv. This conjecture has been proved in \cite{gobbino:ms} by
reducing, via an integral-geometric approach, to the simpler family of
one-dimensional functionals
$$DG_\ep(u)=\frac{1}{\ep}\int_{\re}\arctan\left(
\frac{(u(x+\ep)-u(x))^2}{\ep}\right)\,dx.$$

In this paper we generalize this result. To this end, we consider the
family of functionals
\begin{equation}
	\mathcal{F}_\ep(u)=\int_{\re^n\times\re^n}\varphi_{\ep|\xi|}\left(
	\frac{\left|u(x+\ep\xi)-u(x)\right|}{\ep|\xi|}\right)\eta(\xi)
	\,d\xi\,dx,\label{defn:fep-ndim}
\end{equation}
where $\{\varphi_{\rho}\}_{\rho>0}$ is a family of Borel functions,
and $\eta\in L^{1}(\re^n)$.

Our aim is twofold:
\begin{itemize}
	\item given $\{\varphi_{\rho}\}$, providing estimates for the
	$\Gamma$-limit of $\{\mathcal{F}_{\ep}\}$ in
	terms of $\{\varphi_{\rho}\}$;

	\item  given a functional $\mathcal{F}$ of the form
	(\ref{defn:fdp-gen}), finding $\{\varphi_{\rho}\}$ such that the
	family $\{\mathcal{F}_{\ep}\}$ defined as in (\ref{defn:fep-ndim})
	converges to $\mathcal{F}$.
\end{itemize}

In particular, if $\varphi$ and $\psi$ satisfy the usual assumptions
in order to have lower semicontinuity of $\mathcal{F}$, and $\varphi$
is ``sectionable'' according to Definition \ref{defn:sectionable}
(\eg\ $\varphi(r)=|r|^p$ with $p>1$), then we prove (Theorem\
\ref{thm:main}) that there exists $\{\varphi_{\rho}\}$ such that
the following convergence properties are satisfied:
\begin{enumerate}
	\renewcommand{\labelenumi}{(C\arabic{enumi})\ }

	\item ${\cal F}_\varepsilon(u)\leq\mathcal{F}(u)$ for every
	$u\in L^{1}_{loc}(\re^{n})$;

	\item $\{{\cal F}_\varepsilon(u)\}$ pointwise converges to
	$\mathcal{F}(u)$;

	\item $\mathcal{F}(u)$ is the $\Gamma^{\mbox{-}}$-limit of ${\{\cal
	F}_\varepsilon(u)\}$ in $L^1_{{\rm loc}}(\re^n)$;

	\item if $\sup_{\ep>0}\left\{\mathcal{F}_{\ep}(\uep)+
	\|\uep\|_{\infty}\right\}<+\infty$, then there exist
	$\{\ep_{j}\}\to 0^{+}$ and $u\in \mathit{GSBV}(\re^n)$ such that
	$\left\{u_{\ep_{j}}\right\}\to u$ in $L^1_{{\rm loc}}(\re^n)$.
\end{enumerate}

As in the case of the Mumford-Shah functional, the theory relies
almost completely on the study of the simpler family of
one-dimensional functionals
\begin{equation}
	F_\ep(u)=\int_{\re}\phiep\left(
	\frac{|u(x+\ep)-u(x)|}{\ep}\right)\,dx.
	\label{defn:fep}
\end{equation}

We point out that pointwise estimates like (C1) are one of the main
advantages of this approach.  Thanks to such estimates, the passage
from the one-dimensional to the $n$-dimensional case is a simple
application of Fatou's lemma and standard integral geometric
equalities.

For this reason the finite difference approach is, at the present, the
only approach which has been proved to work also with functionals as
(\ref{defn:fdp-gen}), in the case where $\mathcal{F}(u)$ can be finite
even if $\hn(S_{u})=+\infty$ (this happens \eg\ if $\varphi(r)=r^{2}$
and $\psi(r)=\sqrt{r}$).

This paper is organized as follows: in \S\,\ref{sec:prelim} we give
notations and preliminaries; in \S\,\ref{sec:onedim} we study the
convergence of the functionals defined in
(\ref{defn:fep}); in \S\,\ref{sec:ndim} we consider the general family
(\ref{defn:fep-ndim}) and we prove (C1), (C2), and (C3) under
suitable assumptions on $\{\varphi_{\rho}\}$; in
\S\,\ref{sec:compactness} we consider the compactness property (C4); in
\S\,\ref{sec:approximation} we prove our main approximation result for
the functional $\mathcal{F}(u)$ (Theorem \ref{thm:main}); in
\S\,\ref{sec:examples} we show some simple examples where the theory
developed in this paper applies.

Finally, we would like to thank the referee for carefully reading the
manuscript.

\setcounter{equation}{0}
\section{Preliminaries}
\label{sec:prelim}

In this section we fix notations and we recall basic definitions from
the theory of $\mathit{SBV}$ functions and \gconv.

For all $\alpha\in\re$ the integer part of $\alpha$ is denoted by
$[\alpha]=\sup\{z\in\z\tc z\leq\alpha\}$.  Given $x,y\in\re^{n}$,
their scalar product is denoted by $\langle x,y\rangle$, and the
Euclidean norm of $x$ is denoted by $|x|$.  Given $a,b\in\re$, the
maximum and the minimum of $\{a,b\}$ are denoted by $a\vee b$ and
$a\wedge b$, respectively.  Given $A,B\sse\re^{n}$, we write
$A\subset\subset B$ if the closure of $A$ is compact and contained in
$B$.

The Lebesgue measure and the $(n-1)$-dimensional Hausdorff measure of
a set $B\sse\re^{n}$ are denoted by $|B|$ and $\hn(B)$ respectively.
The restriction of the measure $\hn$ to the set $B$ is denoted by
$\hn\lfloor_{B}$.
We use standard notations for the Banach spaces $L^{p}(\re^{n})$ and
$W^{1,p}(\re^{n})$, and for the metrizable spaces
$L^{p}_{loc}(\re^{n})$.  All the functionals introduced in this paper,
and also all the operations of $\lim$, $\liminf$, $\limsup$, are
intended with range in the extended real line
$\overline{\re}=\re\cup\{+\infty,-\infty\}$.

\subsection{Special functions of bounded variation}
\label{subsec:sbv}
For the general theory of functions with bounded variation we refer to
\cite{eg,ziemer}; here we just recall some definitions and some basic
results.

Let $\Omega\sse\re^{n}$ be an open set, let $u:\Omega\nlra\re$ be a
measurable function, and let $x\in \Omega$. We denote by $u^{+}(x)$
and $u^{-}(x)$, respectively, the upper and lower limit of $u$ at $x$,
defined by
\begin{eqnarray*}
	u^{+}(x) & := & \inf\left\{t\in\re\tc \lim_{\rho\to
	0^{+}}\frac{\left|\left\{y\in\Omega\tc|x-y|<\rho,\
	u(y)>t\right\}\right|}{\rho^n}=0\right\},  \\
	\noalign{\vs}
	u^{-}(x) & := & \sup\left\{t\in\re\tc \lim_{\rho\to
	0^{+}}\frac{\left|\left\{y\in\Omega\tc|x-y|<\rho,\
	u(y)<t\right\}\right|}{\rho^n}=0\right\}.
\end{eqnarray*}

If $u^{+}(x)=u^{-}(x)\in\re$, then $x$ is said to be a \emph{Lebesgue point}
of $u$; in this case, the common value of  $u^{+}(x)$
and $u^{-}(x)$ is called the \emph{approximate limit} of $u$ at the point
$x$, and is denoted by $\mbox{ap\,-}\!\lim_{y\to x}u(y)$. We denote by
$S_{u}$ the \emph{discontinuity set} of $u$, \ie\ the set of all
$x\in\Omega$ which are not Lebesgue points of $u$.

We say that $u$ is a {\em function of bounded variation} in $\Omega$, and
we write $u\in BV(\Omega)$, if $u\in L^{1}(\Omega)$ and its distributional
derivative is a vector-valued measure $Du$ with finite total variation
$|Du|(\Omega)$. We recall that the total variation in $\Omega$ can be
defined also for every measurable function $v:\Omega\to\re$ by the
formula
\begin{equation}
	\left|Dv\right|(\Omega)\ :=\
	\sup\left\{\int_{\Omega}^{}v\,\mathrm{div}\varphi\tc\varphi\in
	C^\infty_{0}(\Omega,\re^n),\ \|\varphi\|_{\infty}\leq 1\right\}.
	\label{defn:du}
\end{equation}

If $u\in BV(\Omega)$, then $S_{u}$ turns out to be countably $(\hn,n-1)$
rectifiable, {\em i.e.}
$$S_{u}=N\cup\bigcup_{i\in\n}K_{i},$$
where $\hn(N)=0$, and each $K_{i}$ is a compact set contained in a
$C^{1}$ hypersurface.

For every $u\in BV(\Omega)$ we have the decomposition
$Du=D^{a}u+D^{s}u$, where $D^{a}u$ is absolutely continuous and
$D^{s}u$ is singular with respect to the Lebesgue measure.  The
density of $D^{a}u$ with respect to the Lebesgue measure is denoted
by $\nabla u$. It turns out that, for almost every $x\in\Omega$, the
vector $\nabla u(x)$ is the \emph{approximate gradient} of $u$, \ie
$$\mbox{ap-}\!\lim_{y\nlra x}\,\frac{u(y)-u(x)-\langle\nabla u(x),y-x
\rangle}{|y-x|}=0.$$

Moreover, we denote the restriction of $D^{s}u$ to $S_{u}$ by
$D^{j}u$, and the restriction of $D^{s}u$ to $\Omega\backslash S_{u}$
by $D^{c}u$.  With these notations we have the following
decomposition:
$$Du=D^{a}u+D^{j}u+D^{c}u.$$

The reader interested in the structure of $D^{a}u$, $D^{j}u$,
$D^{c}u$ is referred to \cite{ambrsc,ambrosiosurvey}.

We say that $u$ is a {\em special function of bounded variation}, and
we write $u\in \mathit{SBV}(\Omega)$, if $u\in BV(\Omega)$ and $D^{c}u=0$.  We
consider also the larger space $\mathit{GSBV}(\Omega)$, which is composed by
all measurable functions $u:\Omega\lra\re$ whose truncations
$u_{k}=(u\wedge k)\vee(-k)$ belong to $\mathit{SBV}(\Omega')$ for every $k>0$,
and every open set $\Omega'\subset\subset \Omega$.

Every $u\in \mathit{GSBV}(\Omega)\cap L^{1}_{loc}(\Omega)$ has an
approximate gradient
$\nabla u(x)$ for a.e. $x\in \Omega$, and a countably $(\hn,n-1)$
rectifiable discontinuity set $S_{u}$.

The spaces $\mathit{SBV}(\Omega)$ and $\mathit{GSBV}(\Omega)$ have been
introduced by De Giorgi
and Ambrosio in \cite{dga}, and have been studied in \cite{ambrarch}.

Given two Borel functions $\varphi:[0,+\infty[\to[0,+\infty]$ and
$\psi:]0,+\infty]\to[0,+\infty]$, we consider the functional
$\mathcal{F}_{\varphi,\psi}:L^{1}_{loc}(\Omega)\to[0,+\infty]$ defined
by
\begin{equation}
	\mathcal{F}_{\varphi,\psi}(u,\Omega)=\left\{ \begin{array}{ll}
	{\displaystyle \int_{\Omega}\varphi(|\nabla u|)\,dx+
	\int_{S_{u}}\psi(|u^{+}-u^{-}|)\,d\hn} & \mbox{if
	}u\in \mathit{GSBV}(\Omega) \\
		\noalign{\vs}
	+\infty & \mbox{otherwise}
	\end{array}\right.
	\label{defn:f-phi-psi}
\end{equation}

This functional is ``isotropic'', in
the sense that it is invariant under rigid motions.  With an abuse of
notation, we denote by $\mathcal{F}_{\varphi,\psi}$ also the
non-isotropic functional where the first integral is replaced by the
integral of $\varphi(\nabla u)$, where $\varphi:\re^n\to[0,+\infty]$.

In \cite{ambrarch} the following semicontinuity result is proved.

\begin{thm}\label{thm:lsc}
	Let $\Omega\sse\re^n$ be an open set.  Let $\varphi:[0,+\infty[ \to
	[0, +\infty]$ be a non-decreasing convex function such that
	\begin{equation}
		\lim_{r\to +\infty} \frac{\varphi(r)}{r} = +\infty,
		\label{hp:phi}
	\end{equation}
	and let $\psi : ]0,+\infty] \to [0,
	+\infty]$ be a non-decreasing concave function such that
	\begin{equation}
		\lim_{r\to 0^{+}}\frac{\psi(r)}{r}= +\infty.
		\label{hp:psi}
	\end{equation}

	Then the functional $\mathcal{F}_{\varphi,\psi}(u,\Omega)$ defined
	in (\ref{defn:f-phi-psi}) is lower semicontinuous in
	$L^{1}_{loc}(\Omega)$.\qed
\end{thm}

\subsection{\gconv}
\label{subsec:gconv}
		Let $X$ be a metric space, let $\{F_{i}\}$ be a sequence of
		functions defined in $X$ with values in $\overline{\re}$.
Let us set
		\begin{eqnarray*}
			\gmliminf{i\to\infty}F_{i}(x) & := &
			\inf\left\{\liminf_{i\to\infty} F_{i}(x_{i})\tc
			\{x_{i}\}\to x\right\}, \\
			\gmlimsup{i\to\infty}F_{i}(x) & := & \inf\left\{
			\limsup_{i\to\infty} F_{i}(x_{i})\tc \{x_{i}\}\to
			x\right\}.
		\end{eqnarray*}

It turns out that $\gmliminf{i\to\infty}F_{i}(x)$ and
$\gmlimsup{i\to\infty}F_{i}(x)$ are lower semicontinuous functions.
Moreover, the ``$\inf$'' in the definitions above are actually
``$\min$''.

If $\gmliminf{i\to\infty} F_{i}(x)=\gmlimsup{i\to\infty}F_{i}(x)=F(x)$
for all $x\in X$, we say that $F$ is the $\Gamma^{\mbox{-}}$-limit of
$\{F_{i}\}$, and we write
$$F(x)=\gmlim{i\nlra\infty}F_{i}(x).$$

This means that for every $x\in X$ the following two conditions are
satisfied:
\begin{enumerate}
	\renewcommand{\labelenumi}{(\roman{enumi})\ }
	\item if
	$\{x_{i}\}\to x$ is any sequence, then
	$\liminf_{i\nlra\infty}F_{i}(x_{i})\geq F(x)$;

	\item there exists a sequence $\{x_{i}\}\to x$ such that
	$F(x)\geq\limsup_{i\nlra\infty}F_{i}(x_{i})$.
\end{enumerate}

The $\Gamma^{\mbox{-}}$-limit, when it exists, is unique, and stable
under subsequences.  The reader interested in variational properties
of \gconv\ is referred to \cite{dmgconv}.

In general, there is no relation between the $\Gamma^{\mbox{-}}$-limit
and the pointwise limit.  However, if $\{F_{i}\}\to F$ uniformly on
compact subsets of $X$, then $F$ is also the
$\Gamma^{\mbox{-}}$-limit of $\{F_{i}\}$.

A special case is when $F_{i}(x)=G(x)$ for every $i\in\n$: in this
case the \mbox{$\Gamma^{\mbox{-}}$-limit} of $\{F_{i}\}$ is the so called
relaxation of $G$, which we denote by $\overline{G}$. We recall that
$\overline{G}$ can also be defined as the supremum of all the lower
semicontinuous functions less or equal than $G$.

Finally, we say that a family $\{F_{\varepsilon}\}_{\ep>0}$ of
functions $\Gamma^{\mbox{-}}$-converges to $F$ as $\ep\nlra 0^+$, if
$\{F_{\varepsilon_{i}}\}$ $\Gamma^{\mbox{-}}$-converges to $F$ for
every sequence $\{\ep_{i}\}\nlra 0^{+}$.

\setcounter{equation}{0}
\section{The One-Dimensional Functionals $F_{\ep}$}
\label{sec:onedim}
In this section we consider a family $\{\phiep\}_{\ep>0}$ of Borel
functions $\phiep:[0,+\infty[\nlra [0,+\infty[$, and we study the
convergence of the family of functionals
\begin{equation}
	\fep(u,\Omega):=\int_{\Omega}\phiep\left(
	\frac{|u(x+\ep)-u(x)|}{\ep}\right)\,dx,
	\label{defn:fep-omega}
\end{equation}
defined for every $\ep>0$, $u\in L^{1}_{loc}(\re)$, and every
measurable set $\Omega\sse\re$, with values in $\re\cup\{+\infty\}$.
When $\Omega=\re$, then we simply write $\fep(u)$ instead of
$\fep(u,\re)$.

\subsection{Statement of the results}
We state here all the results which will be proved in this section.

The first one provides an estimate from below for the
$\Gamma^{\mbox{-}}$-limit of $\{F_{\ep}\}$.

\begin{thm} \label{thm:liminf}
Let $\{\phiep\}_{\ep>0}$ be a family of Borel functions such that
\begin{description}
	\item[(\textnormal{li1})] $\phiep$ is
	continuous and non-decreasing for every $\ep>0$;

	\item[(\textnormal{li2})] for each $\ep>0$, the function $\phiep$ is
	either convex, or concave, or ``convex-concave'', \ie\ there
	exists $\overline{r}_{\ep}>0$ such that $\phiep$ is
	convex in $[0, \overline{r}_{\ep}]$ and concave in
$[\overline{r}_{\ep},
	+\infty[$.
\end{description}

Let us define $\varphi_{\star}:[0,+\infty[\to[0,+\infty]$ and
$\psi_{\star}:]0,+\infty]\to[0,+\infty]$ by
\begin{equation}
	\varphi_{\star}(r):=\gmliminf{\ep\to 0^{+}}\,\phiep(r),\hspace{3em}
	\psi_{\star}(r):=\gmliminf{\ep\to 0^{+}}\,\ep\phiep
	\left( \frac{r}{\ep} \right).
	\label{defn:phi-psi}
\end{equation}

Then for every $u\in L^{1}_{loc}(\re)$ we have that
$$\gmliminf{\ep \to 0^{+}}\,F_{\ep} (u,\re)\ \geq\
\overline{\mathcal{F}_{\varphi_{\star},\psi_{\star}}}(u,\re),$$
where $\overline{\mathcal{F}_{\varphi_{\star},\psi_{\star}}}$ is the
relaxation of the functional
$\mathcal{F}_{\varphi_{\star},\psi_{\star}}$ defined as in
(\ref{defn:f-phi-psi}).
\end{thm}

The following result provides a pointwise estimate from above for
$F_{\ep}(u,\re)$.

\begin{thm}\label{thm:estimate}
	Let $\{\phiep\}_{\ep>0}$ be a family of Borel functions such that
	\begin{description}
		\item[\textnormal{(li1)}] $\phiep$ is continuous and
		non-decreasing for every $\ep>0$;

		\item[\textnormal{(Est)}] there exist a convex function
		$\varphi^\star:[0,+\infty[\to[0,+\infty]$, and a concave
		function $\psi^\star:]0,+\infty]\to[0,+\infty]$ such that
		$$\phiep(A+S)\ \leq\ \varphi^\star(A)+
		\frac{1}{\ep}\,\psi^\star(\ep S),$$
		for every $\ep>0$, $A\geq 0$, $S>0$.
	\end{description}

	Then $\fep(u,\re)\leq\overline{\mathcal{F}_{\varphi^{\star},
	\psi^{\star}}}(u,\re)$ for every $\ep>0$, and every $u\in
	L^{1}_{loc}(\re)$.
\end{thm}

In many cases, the pointwise limit and the $\Gamma^{\mbox{-}}$-limit
of $\{F_{\ep}\}$ are uniquely determined by Theorem \ref{thm:liminf}
and Theorem \ref{thm:estimate}, as in the following situation.

\begin{cor}
	\label{cor:liminf-estimate}
	Let us assume that the family $\{\phiep\}$ satisfies assumptions
	(li1), (li2), (Est), and that
	$\varphi_{\star}=\varphi^{\star}=:\varphi$ and
	$\psi_{\star}=\psi^{\star}=:\psi$.

	Then
	\begin{enumerate}
		\renewcommand{\labelenumi}{(\roman{enumi})\ }
		\item
		$\fep(u,\re)\leq\overline{\mathcal{F}_{\varphi,\psi}}(u,\re)$
		for every $u\in L^{1}_{loc}(\re)$, and every $\ep>0$;

		\item $\{\fep(u,\re)\}$ pointwise converges to
		$\overline{\mathcal{F}_{\varphi,\psi}}(u,\re)$;

		\item $\overline{\mathcal{F}_{\varphi,\psi}}(u,\re)$ is the
		$\Gamma^{\mbox{-}}$-limit of $\{F_{\ep}(u,\re)\}$ in
		$L^{1}_{loc}(\re)$.
	\end{enumerate}
\end{cor}

\begin{rmk}
	\begin{em}
		All the results stated above can be generalized
word-by-word to the
		vector valued case $u\in L^{1}_{loc}(\re;\re^k)$.
	\end{em}
\end{rmk}

\subsection{Estimates from below}
In this subsection we prove Theorem \ref{thm:liminf}. The strategy of
the proof follows the argument used in \cite{gobbino:ms} in the case
of the Mumford-Shah functional.

In order to avoid a cumbersome notation, we extend to $[0,+\infty]$
the function $\psi_{\star}$, defined in (\ref{defn:phi-psi}), by setting
$\psi_{\star}(0)=0$.  Moreover, for every $\alpha\geq 0$, $\beta>0$,
we define
\begin{equation}
	\lambda(\alpha, \beta)\ :=\ \min\left\{\beta \varphi_{\star}
	\left(\frac{\alpha - l}{\beta}\right) + \psi_{\star} (l)\tc 0\leq
	l\leq \alpha\right\}.
	\label{defn:lambda}
\end{equation}

Since $\varphi_{\star}$ and the extension of $\psi_{\star}$ are lower
semicontinuous
on $[0,+\infty[$, it follows that $\lambda$ is well defined and lower
semicontinuous. Moreover, since we can always set
$l=0$ in (\ref{defn:lambda}), it turns out that
\begin{equation}
	\lambda(\alpha, \beta)\ \leq\ \beta\,\varphi_{\star}
	\left(\frac {\alpha}{\beta}\right).
	\label{estim:lambda}
\end{equation}

Now we state and prove three technical lemmata.  The first one is a
``discretization'' of Theorem \ref{thm:liminf}.
\begin{lemma}\label{lemma:anal1}
Let $\{\varphi_{\ep}\}_{\ep >0}$ be a family of Borel functions
satisfying (li1), (li2), and let $\alpha \geq 0$ and $\beta >0$.
Then for every $\ep \in ]0, \beta]$, there exists
\begin{equation}
	\label{defn:theta}
	\Theta (\ep, \alpha, \beta):= \min \left\{ \sum_{i=1}^{N_{\ep}} \ep
	\varphi_{\ep} \left(\frac{|x_{i}|}{\ep}\right) \tc
	\sum_{i=1}^{N_{\ep}} |x_{i}| \geq \alpha, \ N_{\ep}
	= \left[ \frac{\beta}{ \ep} \right] \right\}.
\end{equation}

Moreover,
\begin{equation}
	\liminf_{\ep \to 0^{+}} \,\Theta (\ep, \alpha, \beta)\ \geq\
	\lambda(\alpha, \beta),
	\label{tslemma:anal1}
\end{equation}
where $\lambda$ is the function defined in (\ref{defn:lambda}).
\end{lemma}

{\sc Proof.}
The minimum problem (\ref{defn:theta}) has at least one solution,
since by (li1) we can restrict to the compact set
$$\left\{ (x_{1}, x_{2}, \ldots, x_{N_{\ep}}) \in [0,
\alpha]^{N_{\ep}} \tc \sum_{i=1}^{N_{\ep}} x_{i}= \alpha \right\}.$$

Therefore the function $\Theta(\ep,\alpha,\beta)$ is well defined.

Now let $\{\ep _{n}\} \to 0^{+}$ be a sequence such that
$$\liminf_{\ep \to 0^{+}}\, \Theta (\ep, \alpha,
\beta)\ =\ \lim_{n\to \infty}\,\Theta (\ep_{n}, \alpha, \beta),$$
and, for all $\ep_{n}$, let $x_{n,1}\geq x_{n,2}\geq \ldots \geq
x_{n,N_{\ep_{n}}}$ be a minimizer for (\ref{defn:theta}).

Since by (li2) the function $r \mapsto \ep_{n} \varphi_{\ep_{n}}
\left(r/\ep_{n}\right)$ is convex in $[0,
\ep_{n}\overline{r}_{\ep_{n}}]$ and concave in
$[\ep_{n}\overline{r}_{\ep_{n}}, +\infty[$ (with obvious modifications
if $\phiep$ is always convex or always concave), it follows that only
$x_{n,1}$ can be greater than $\ep_{n}\overline{r}_{\ep_{n}}$, and all
the $x_{n,i}$'s in the convexity zone are equal (this is true if
$\varphi_{\ep}$ is strictly convex in $[0,\overline{r}_{\ep}]$;
however, if $\varphi_{\ep}$ has a flat zone in
$[0,\overline{r}_{\ep}]$, then there exists at least one minimizer
with the property that all the $x_{n,i}$'s in the convexity zone are
equal, and so we can work with this minimizer without loss of
generality).  Therefore, there are only two possibilities:
\begin{description}
	\item [\textnormal{(P1)}] $x_{n,1} =\ldots= x_{n, N_{\ep_{n}}}=
	\alpha/N_{\ep_{n}}$, and in this case
	\begin{equation}
		\Theta (\ep_{n}, \alpha, \beta)\ =\ \ep_{n} N_{\ep_{n}}
		\varphi_{ \ep_{n}} \left(\frac{\alpha}{ \ep_{n}
		N_{\ep_{n}}}\right);
		\label{formula:theta1}
	\end{equation}

	\item [\textnormal{(P2)}] $x_{n,1}\geq
\ep_{n}\overline{r}_{\ep_{n}}$ and
	$x_{n,2}=\ldots = x_{n,N_{\ep_{n}}}= (\alpha -
	x_{n,1})/(N_{\ep_{n}}-1)$.  In this case
	\begin{equation}
		\Theta (\ep_{n}, \alpha, \beta)\ =\ \ep_{n} \varphi_{ \ep_{n}}
		\left(\frac{x_{n,1}}{ \ep_{n}}\right) + \ep_{n} (
		N_{\ep_{n}}-1)\, \varphi_{ \ep_{n}} \left(\frac{\alpha -
		x_{n,1}}{ \ep_{n} (N_{\ep_{n}}-1)}\right).
		\label{formula:theta2}
	\end{equation}
\end{description}

Up to subsequences, we can suppose that either (P1) or (P2) holds true for all
$n\in\n$. In the first case, observing that $\{\ep_{n}
N_{\ep_{n}}\}\to\beta$ and using the definition of $\varphi_{\star}$,
passing to the limit in (\ref{formula:theta1}) we have that
$$\liminf_{n\to \infty}\, \Theta (\ep_{n}, \alpha, \beta)\ \geq\ \beta\,
\varphi_{\star}\left( \frac {\alpha}{\beta} \right)\ \geq\
\lambda(\alpha,\beta).$$

In the second case, up to subsequences, we can assume that there
exists
$$l= \lim_{n\to \infty} x_{n,1} \in [0, \alpha].$$

By the definition of $\varphi_{\star}$, $\psi_{\star}$, and $\lambda$,
passing to the limit in (\ref{formula:theta2}) we obtain that
$$\liminf_{n\to \infty}\, \Theta (\ep_{n}, \alpha, \beta)\ \geq\
\psi_{\star}(l) + \beta \varphi_{\star}\left(\frac {\alpha - l}{\beta}
\right)\ \geq\ \lambda (\alpha, \beta).$$

In both cases, inequality (\ref{tslemma:anal1}) is proved.
\qed

The second lemma is a ``localization'' of Theorem \ref{thm:liminf}.
\begin{lemma}\label{lemma:local}
Let $I= [a,b]$ be an interval, let $\{u_{\ep}\}\sse
L^{1}_{loc}(\re)$, and let $u \in L^{1}_{loc}(\re)$.  Let us assume
that
\begin{enumerate}
	\renewcommand{\labelenumi}{(\roman{enumi})\ }
	\item $u_{\ep} \to u$ in $L^{1}_{loc}(\re)$;

	\item $a$ and $b$ are Lebesgue points of $u$.
\end{enumerate}
Then
\begin{equation}
	\liminf_{\ep \to 0^{+}} F_{\ep} (u_{\ep}, I)\ \geq\ \lambda (|u(b) -
	u(a) |,b-a).
	\label{tslemma:local}
\end{equation}
\end{lemma}
\prf\ Let $\{\ep_{n}\}\to 0^{+}$ be a sequence such that
$$\liminf_{\ep \to 0^{+}} F_{\ep} (u_{\ep}, I)\ =\
\lim_{n\to \infty} F_{\ep_{n}} (u_{\ep_{n}}, I).$$
Up to subsequences, we can assume that
$$u_{\ep_{n}}(x)\lra u(x)\hspace{2em}\mbox{for a.e. }x\in I.$$

Now, let us set $J:=|u(b) - u(a)|$.  If $J=0$, then there is nothing
to prove.  Otherwise, let us fix $\eta\in ]0,J]$, let us set
$N_{\ep_{n}} = \left[ (b-a)/\ep_{n} \right]$, and let us define
$$C_{n} = \left\{ x \in [a, a + \ep_{n}] :
\sum_{k=1}^{N_{\ep_{n}}} |u_{\ep_{n}} (x + k\ep_{n}) - u_{\ep_{n}}(x +
(k-1)\ep_{n})| \geq J- \eta \right\}.$$

Using assumption (ii), it can be proved (for the technical details see
Step 2 in the proof of Lemma 3.2 in \cite{gobbino:ms}) that
\begin{equation} \label{cn}
\lim_{n\to \infty} \frac{|C_{n}|}{\ep_{n}}\ =\ 1.
\end{equation}

By the definition of $C_{n}$ we obtain that
\begin{eqnarray*}
	F_{\ep_{n}}(u_{\ep_{n}}, I) & \geq &
	F_{\ep_{n}}(u_{\ep_{n}}, [a, a+ \ep_{n} N_{\ep_{n}}]) \\
	 \noalign{\vs}
	& = & \int_{a}^{a+\ep_{n}N_{\ep_{n}}}
	 \varphi_{\ep_{n}} \left( \frac{|u_{\ep_{n}}(x+\ep_{n}) -
	 u_{\ep_{n}} (x)|}{\ep_{n}} \right) dx \\
	\noalign{\vs}
	& = & \frac{1}{\ep_{n}}
	\int_{a}^{a +\ep_{n}} \sum_{k=1}^{N_{\ep_{n}}}
	\ep_{n} \varphi_{\ep_{n}}\left(
	\frac{|u_{\ep_{n}}(x+ k \ep_{n}) -
	u_{\ep_{n}} (x+ (k-1) \ep_{n})|}{\ep_{n}}
	\right) dx \\
	\noalign{\vs}
	 & \geq & \frac{|C_{n}|}{\ep_{n}}\, \Theta (\ep_{n},
	J-\eta,b-a)
\end{eqnarray*}
where $\Theta$ is the function defined in (\ref{defn:theta}). Applying
Lemma \ref{lemma:anal1} with $\alpha = J-\eta$ and $\beta=b-a$, and using
(\ref{cn}), we conclude that
$$\liminf_{n\to \infty}\,F_{\ep_{n}} (u_{\ep_{n}}, I)\ \geq\
\liminf_{n\to \infty}\,\,\frac{|C_{n}|}{\ep_{n}}\, \Theta
(\ep_{n}, J-\eta,b-a)\ \geq\ \lambda (J-\eta,b-a).$$

Since $\lambda$ is lower semicontinuous, and $\eta$ is arbitrary,
(\ref{tslemma:local}) is proved.
\qed

The third Lemma states a general property of $L^{p}(\re)$ spaces (for
a proof, see Lemma 3.3 in \cite{gobbino:ms}).
\begin{lemma} \label{lemma:anal3}
Let $u\in L^{\infty}(\re)$. Then there exists $a\in \re$ such that
\begin{enumerate}
	\renewcommand{\labelenumi}{(\roman{enumi})\ }
	\item $a + q$ is a Lebesgue point of $u$ for every rational number $q$;

	\item every sequence $\{u_{n}\}\subset L^{\infty}(\re)$ which
	satisfies the following two conditions
	\begin{itemize}
	 	\item  $u_{n}(a + \frac{z}{n})= u(a + \frac{z}{n})$ for all
	 	$z\in\z$,

	 	\item  if $x\in [a+\frac{z}{n}, a+\frac{(z+1)}{n}]$, then
		$u_{n}(x)$ belongs to the interval with endpoints
		$u(a+\frac{z}{n})$
	 	and $u(a+\frac{(z+1)}{n}),$
	 \end{itemize}
has a subsequence converging to $u$ in $L^{1}_{loc}(\re).$
\qed
\end{enumerate}
\end{lemma}
\textsc{Proof of Theorem \ref{thm:liminf}.}
Let us set for simplicity
$F_{\star}:=\mathcal{F}_{\varphi_{\star},\psi_{\star}}$. We have to show that
$$\liminf_{n\to \infty} F_{\ep_{n}} (u_{n}) \geq \overline{F_{\star}}(u)$$
for every $u\in L^{1}_{loc}(\re)$, every $\{\ep_{n}\}\to 0^{+}$, and
every sequence $\{u_{n}\} \to u$ in $L^{1}_{loc}(\re)$.  Let us begin
with the case where $u\in
L^{\infty}(\re)$, $\{u_{n}\}\sse L^{\infty}(\re)$, and
$\|u_{n}\|_{\infty}\leq\|u\|_{\infty}$.

Our strategy is to construct a sequence $\{v_{j}\}\sse \mathit{GSBV}(\re) \cap
L^{1}_{loc}(\re)$ such that
\begin{equation}
	\label{punto1}
	\{v_{j}\}\to u\hspace{2em}\hbox{in}\; L^{1}_{loc}(\re),
\end{equation}
\begin{equation}
	\label{punto2}
	\liminf_{n \to \infty} F_{\ep_{n}} (u_{n}) \geq F_{\star}(v_{j})
	\hspace{2em}\forall j \in\n.
\end{equation}

Let us assume that $a\in [0,1]$ satisfies conditions (i) and (ii) of
Lemma \ref{lemma:anal3}.  For all $z\in\z$ and $j \in\n$, let
$\beta=1/j$, let
$I_{j}^{z}$ be the interval $[a + z\beta, a + (z+1)\beta]$, let $J_{j}^{z}$
be the increment $|u(a + (z+1)\beta) - u(a + z\beta)|$, and let $l_{j}^{z}$
be such that
$$\lambda(J_{j}^{z},\beta)\ =\ \beta \varphi_{\star} \left(\frac{J_{j}^{z}
-l_{j}^{z}}{\beta}\right) + \psi_{\star} (l_{j}^{z}).$$

Now we define $v_{j}$ on every interval $I_{j}^{z}$ as the piecewise
affine function which
\begin{itemize}
	\item coincides with $u$ at the endpoints of $I_{j}^{z}$;

	\item has constant (approximate) gradient in the interval, with
	$|\nabla v_{j}(x)|=(J_{j}^{z} -l_{j}^{z})/\beta$;

	\item has a jump of height $l_{j}^{z}$ in the medium point of the
	interval (of course if $l_{j}^{z}=0$, then no jump point is
	necessary).
\end{itemize}

It is easy to check that the functions $v_{j}$ satisfy both
assumptions of (ii) of Lemma \ref{lemma:anal3}; hence, up to
subsequences, (\ref{punto1}) holds.  Moreover, $v_{j}$ belongs to
$\mathit{GSBV}(\re)\cap L^{1}_{loc}(\re)$, and for all $z\in\z$, $j\in\n$,
we have that
$$F_{\star}(v_{j}, I_{j}^{z})\ =\ \lambda(J_{j}^{z}, \beta),$$
hence, by Lemma \ref{lemma:local} applied in the interval $I_{j}^{z}$,
$$\liminf_{n\to \infty} F_{\ep_{n}} (u_{n}, I_{j}^{z})\ \geq\
\lambda(J_{j}^{z}, \beta)\ =\
F_{\star}(v_{j}, I_{j}^{z}).$$

Summing over all $z\in\z$ and using Fatou's Lemma for series, it
follows that (\ref{punto2}) holds true, and this completes the proof
in the case $u\in L^{\infty}(\re)$.

In the general case $u\in L^{1}_{loc}(\re)$, let us denote by $T_{k}$
the truncation operator $T_{k}v=(v\vee k)\wedge k$. For every $k>0$ we
have that $\{T_{k}u_{n}\}\to T_{k}u$ as $n\to +\infty$.
Moreover, since $\varphi_{\ep}$ is non-decreasing we have that
$F_{\ep_{n}}(u_{n},\re)\geq F_{\ep_{n}}(T_{k}u_{n},\re)$ so that
\begin{equation}
	\liminf_{n\to \infty} F_{\ep_{n}} (u_{n},\re)\ \geq\
	\liminf_{n\to \infty} F_{\ep_{n}} (T_{k}u_{n},\re)\ \geq\
	F_{\star}(T_{k}u,\re),
	\label{eqn:tk}
\end{equation}
where the last inequality follows from the $L^{\infty}$ case proved
above.

Since $\{T_{k}u\}\to u$ in $L^{1}_{loc}(\re)$ as $k\to +\infty$, then
the conclusion follows letting $k\to +\infty$, due to (\ref{eqn:tk})
and the lower semi-continuity of $F_{\star}$.
\qed

\subsection{Pointwise estimates}
In this subsection we prove Theorem \ref{thm:estimate} and Corollary
\ref{cor:liminf-estimate}.\vs

\noindent\textsc{Proof of Theorem \ref{thm:estimate}.} Since $F_{\ep}$
is lower semicontinuous in $L^{1}_{loc}(\re)$ (by Fatou's Lemma), then
it is enough to prove that
$F_{\ep}(u)\leq\mathcal{F}_{\varphi^\star,\psi^\star}(u)$.  To this
end, we can of course assume that
$\mathcal{F}_{\varphi^\star,\psi^\star}(u)<+\infty$, hence $u\in
\mathit{SBV}(J)$ for every $J\subset\subset\re$.  In this case let us
set, for every $x\in\re$,
\begin{eqnarray*}
	A_{\ep}(x) & := & {\displaystyle \left|D^{a}u\right|([x,x+\ep])\ =\
	\int_{0}^{\ep}\left|\nabla u(x+\tau)\right|\,d\tau;}  \\
	\noalign{\vs}
	S_{\ep}(x) & := & {\displaystyle \left|D^{j}u\right|([x,x+\ep])\ =\
	\int_{x}^{x+\ep}\left|u^{+}(\tau)-u^{-}(\tau)\right|\,d
	\mathcal{H}^{0}\lfloor_{S_{u}}(\tau).}
\end{eqnarray*}

Since $|u(x+\ep)-u(x)|\leq A_{\ep}(x)+S_{\ep}(x)$ for a.e. $x\in\re$, by
(Est) we have that
\begin{equation}
	\phiep\left(\frac{|u(x+\ep)-u(x)|}{\ep}\right)\ \leq\
	\varphi^\star\left(\frac{A_{\ep}(x)}{\ep}\right)
	+\frac{1}{\ep}\,\psi^\star\left(S_{\ep}(x)\right).
	\label{phiep-split}
\end{equation}

Now let us estimate separately the integral of the two summands.
Since $\varphi^\star$ is convex, by Jensen's inequality we have
that
\begin{eqnarray}
	\int_{\re}\varphi^\star\left(\frac{A_{\ep}(x)}{\ep}\right)dx & = &
	\int_{\re}^{}\varphi^\star\left(\frac{1}{\ep}
	\int_{0}^{\ep}\left|\nabla u(x+\tau)\right|d\tau\right)dx
	\nonumber  \\
	\noalign{\vs}
	 & \leq &
	 \int_{\re}^{}\frac{1}{\ep}\int_{0}^{\ep}\varphi^\star\left(
	 \left|\nabla u(x+\tau)\right|\right)d\tau\, dx
	 \nonumber \\
	\noalign{\vs}
	 & = &
\frac{1}{\ep}\int_{0}^{\ep}d\tau\int_{\re}^{}\varphi^\star\left(
	 \left|\nabla u(x+\tau)\right|\right)dx
	 \nonumber \\
	\noalign{\vs}
	 & = & \int_{\re}\varphi^\star\left(\left|\nabla
	 u(x)\right|\right)dx
	\label{estim:phi}
\end{eqnarray}

Since $\psi^{\star}$ is subadditive, then
\begin{eqnarray}
	\frac{1}{\ep}\int_{\re}^{}\psi^{\star}\left(S_{\ep}(x)\right)dx & = &
	\frac{1}{\ep}\int_{\re}\psi^{\star}\left(
	\int_{x}^{x+\ep}\left|u^{+}(\tau)-u^{-}(\tau)\right|\,d
	\mathcal{H}^{0}\lfloor_{S_{u}}(\tau)\right)dx \nonumber   \\
	 \noalign{\vs}
	 & \leq & \frac{1}{\ep}\int_{\re}\left(
	 \int_{x}^{x+\ep}\psi^{\star}\left(\left|
	 u^{+}(\tau)-u^{-}(\tau)\right|\right)\,d
	 \mathcal{H}^{0}\lfloor_{S_{u}}(\tau)\right)dx \nonumber \\
	\noalign{\vs}
	 & = & \frac{1}{\ep}\int_{\re}\psi^{\star}
	 \left(\left|u^{+}(\tau)-u^{-}(\tau)\right|\right)\,d
	\mathcal{H}^{0}\lfloor_{S_{u}}(\tau)
	\int_{\tau-\ep}^\tau dx \nonumber \\
	\noalign{\vs}
	 & = & \int_{S_{u}}\psi^{\star}
	 \left(\left|u^{+}(\tau)-u^{-}(\tau)\right|\right)\,d
	\mathcal{H}^{0}(\tau).
	\label{estim:psi}
\end{eqnarray}

By (\ref{phiep-split}), (\ref{estim:phi}), and (\ref{estim:psi}),
thesis is proved.
\qed\vs

\noindent\textsc{Proof of Corollary \ref{cor:liminf-estimate}.} The
family $\{\phiep\}$ satisfies assumptions (li1) and (Est) of Theorem
\ref{thm:estimate} with $\varphi=\varphi^\star$ and $\psi=\psi^\star$.
This proves statement (i), and in particular
\begin{equation}
	\limsup_{\ep\to 0^{+}}\,F_{\ep}(u)\ \leq\
	\overline{\mathcal{F}_{\varphi,\psi}}(u),\hspace{2em}\foralll u\in
	L^{1}_{loc}(\re).
	\label{dis:limsup}
\end{equation}

Moreover, $\{\phiep\}$ satisfies assumptions (li1) and (li2) of Theorem
\ref{thm:liminf}. Since $\varphi=\varphi_\star$ and
$\psi=\psi_\star$, it follows that
\begin{equation}
  	\gmliminf{\ep\to 0^{+}}F_{\ep}(u)\ \geq\
	\overline{\mathcal{F}_{\varphi,\psi}}(u),\hspace{2em}\foralll u\in
	L^{1}_{loc}(\re).
  	\label{dis:liminf}
\end{equation}

By (\ref{dis:limsup}) and (\ref{dis:liminf}), statements (ii) and (iii)
follow.\qed

\setcounter{equation}{0}
\section{The general family $\mathcal{F}_{\ep}$}
\label{sec:ndim}
In this section we consider a family $\{\phiep\}_{\ep>0}$ of Borel
functions as in \S\ \ref{sec:onedim}, and a non-negative function
$\eta\in L^{1}(\re^n)$. We study the convergence of the family of functionals
\begin{equation}
	\mathcal{F}_\ep(u,\re^n)=\int_{\re^n\times\re^n}\varphi_{\ep|\xi|}
	\left(\frac{\left|u(x+\ep\xi)-u(x)\right|}{\ep|\xi|}
	\right)\eta(\xi)\,d\xi\,dx,
	\label{defn:fep-ndim-bis}
\end{equation}
defined for every $\ep>0$, and every $u\in L^{1}_{loc}(\re^n)$.

The advantage of $\{\mathcal{F}_{\ep}\}$ with respect to the family
$\{\fep\}$ introduced in \S\ \ref{sec:onedim} is twofold:
\begin{itemize}
	\item  it can be defined in every space dimension;

	\item  it fulfills the compactness properties stated in \S\
	\ref{sec:compactness} (the family $\{\fep\}$, on the contrary,
	satisfies no compactness properties).
\end{itemize}

However, the results of \S\ \ref{sec:onedim} are a fundamental tool in
the study of the convergence of $\{\mathcal{F}_{\ep}\}$, due to
integral geometric techniques. To this end, we introduce the
functionals
\begin{equation}
	\label{defn:fepxi}
	F_{\ep, \xi} (u,\re^n) = \int_{\re^n} \varphi_{\ep |\xi|} \left(
	\frac{|u(x + \ep \xi) - u(x) |}{\ep |\xi|} \right) d x,
\end{equation}
defined for every $\ep>0$, $\xi\in\re^n\setminus\{0\}$, $u\in L^1_{{\rm
loc}}(\re^n)$. With this notation
\begin{equation}
	\mathcal{F}_{\ep}(u,\re^n)=\int_{\re^{n}}F_{\ep,\xi}(u,\re^n)
	\,\eta(\xi)\,d\xi
	\hspace{2em}\foralll u\in L^1_{{\rm loc}}(\re^{n}).
	\label{rel:fep-fepxi}
\end{equation}

Now let $\xi\in\re^n\setminus\{0\}$, and let
$\langle\xi\rangle^{\bot}=\{z\in\re^{n}\tc\langle\xi,z\rangle=0\}$ be
the orthogonal space to $\xi$.  For every
$y\in\langle\xi\rangle^{\bot}$ let us consider the function
$u_{\xi,y}:\re\nlra\re$ defined by
\begin{equation}
	u_{\xi,y}(t)=u\left(y+t\frac{\xi}{|\xi|}\right),
	\hspace{2em}\foralll t\in\re.
	\label{defn:uxiy}
\end{equation}

With the substitution $x=y+t\xi/|\xi|$, relation (\ref{defn:fepxi}) can be
rewritten as
\begin{eqnarray}
	F_{\ep,\xi}(u,\re^n) & = & \int_{\langle\xi\rangle^{\bot}}\int_{\re}
	\varphi_{\ep|\xi|}\left(\frac{\left|
	u(y+t\xi/|\xi|+\ep\xi)-u(y+t\xi/|\xi|)\right|}{\ep|\xi|}\right)dt\,dy
	\nonumber  \\
	\noalign{\vs}
	 & = & \int_{\langle\xi\rangle^{\bot}}\int_{\re}
	\varphi_{\ep|\xi|}\left(\frac{\left|
	u_{\xi,y}(t+\ep|\xi|)-u_{\xi,y}(t)\right|}{\ep|\xi|}\right)dt\,dy
	\nonumber  \\
	\noalign{\vs}
	 & = & \int_{\langle\xi\rangle^{\bot}}F_{\ep|\xi|}(u_{\xi,y},\re)\,dy
	\label{rel:fepxi-fep-ndim}
\end{eqnarray}
where $\{\fep\}$ is the family defined in (\ref{defn:fep-omega}).

Thanks to (\ref{rel:fep-fepxi}) and (\ref{rel:fepxi-fep-ndim}), the
functional $\mathcal{F}_{\ep}(u,\re^n)$ can be written in terms of the
one-dimensional sections of $u$.

We now need the following result about one-dimensional
sections of $\mathit{GSBV}$ functions.
\begin{lemma}
	Let $\varphi$ and $\psi$ be as in the lower semicontinuity theorem
	\ref{thm:lsc}.
	\begin{enumerate}
	\renewcommand{\labelenumi}{(\roman{enumi})\ }
		\item Let $u\in \mathit{GSBV}(\re^n)$.  Then for all
$\,\xi\in\re^n$ we
		have that $u_{\xi,y}\in \mathit{GSBV}(\re)$ for a.e.
		$y\in\langle\xi\rangle^{\bot}$, and moreover
	\begin{eqnarray}
		\nabla u_{\xi,y}(t) & = & \langle \nabla u(y+t\xi/|\xi|),
		\xi/|\xi|\rangle,\hspace{2em}{\rm for\ a.e.\ }t\in\re;
		\label{section1}\\
		\noalign{\vs}
		S_{u_{\xi,y}} & = & \{t\in\re\tc y+t\xi/|\xi|\in S_{u}\};
		\label{section2}
	\end{eqnarray}
	\begin{equation}
		u_{\xi,y}^{+}(t)=u^{+}\left(y+t\xi/|\xi|\right),
		\hspace{1em}
		u_{\xi,y}^{-}(t)=u^{-}\left(y+t\xi/|\xi|\right)
		\hspace{1em}\foralll t\in\re.
		\label{section3}
	\end{equation}

		\item Vice-versa: let $u\in L^1_{loc}(\re^n)$, and let
		$\{\xi_1,\ldots,\xi_n\}\sse\re^n$ be a set of linearly
		independent vectors.  If
	\begin{equation}
		\int_{\langle\xi_i \rangle^{\bot}}
		\mathcal{F}_{\varphi,\psi}(u_{\xi_i,y},\re)\,dy<+\infty
	\label{section4}
	\end{equation}
	for all $i\in\{1,\ldots,n\}$, then $u\in \mathit{GSBV}(\re^n)$.

		\item If $\eta\in L^{1}(\re^n)$ is a non-negative non-zero
		function, then for every $u\in L^{1}_{loc}(\re^n)$ we have
		that
		\begin{equation}
			\int_{\re^n}\left(\int_{\langle\xi\rangle^{\bot}}
	\mathcal{F}_{\varphi,\psi}(u_{\xi,y},\re)\,dy\right)\eta(\xi)\,d\xi
			\ =\ \omega\,\mathcal{F}_{S\varphi,\psi}(u,\re^n),
			\label{rel:fsec-f}
		\end{equation}
		where $\omega:=\|\eta\|_{L^{1}(\re)}$, and
		\begin{equation}
			\left(S\varphi\right)(z):=
			\frac{1}{\omega}\int_{\re^n}\varphi\left(\left|\langle
			z,\xi/|\xi|\rangle\right|\right)\eta(\xi)\,d\xi.
			\label{defn:hatphi}
		\end{equation}
	\end{enumerate}
\end{lemma}
\prf
Statements (i) and (ii) follow from \cite[Theorem 3.3]{ambrsc}.

In order to prove (iii), let us assume first that $u\in
\mathit{GSBV}(\re^n)\cap L^{1}_{loc}(\re^n)$.  In this case by
(\ref{section1}), (\ref{section2}), and (\ref{section3}) we have that
\begin{eqnarray*}
	\lefteqn{\hspace{-1em}\int_{\langle\xi\rangle^{\bot}}
	\mathcal{F}_{\varphi,\psi}(u_{\xi,y},\re)\,dy\
	\ =\ \int_{\langle\xi\rangle^{\bot}}\int_{\re}
	\varphi\left(\left|\langle\nabla u
	\left(y+t\frac{\xi}{|\xi|}\right),
	\frac{\xi}{|\xi|}\rangle\right|\right)dt\,dy+} \\
		\noalign{\vs}
	 \hspace{1em}&  &
	 +\int_{\langle\xi\rangle^{\bot}}\int_{S_{u_{\xi,y}}}
	 \psi\left(\left|u_{\xi,y}^{+}(t)-u_{\xi,y}^-(t)
	 \right|\right)d\mathcal{H}^0(t)\,dy  \\
		\noalign{\vs}
	 & = & \int_{\re^n}\varphi\left(\left|\langle
	 \nabla u(x),\xi/|\xi|\rangle\right|\right)dx +\int_{S_{u}}
	 \psi\left(\left|u^{+}(y)-u^-(y)
	 \right|\right)d\hn(y),
\end{eqnarray*}
where the last equality follows from the substitution
$x=y+t\xi/|\xi|$ for the first summand, and from
\cite[Theorem 3.2.26]{federer} for the second summand. Multiplying this
equality by $\eta(\xi)$, and integrating in $\xi$ over $\re^n$, we
prove (\ref{rel:fsec-f}) in this case.

If $u\in L^{1}_{loc}(\re^n)\setminus \mathit{GSBV}(\re^n)$, then necessarily
$$\int_{\langle\xi\rangle^{\bot}}
\mathcal{F}_{\varphi,\psi}(u_{\xi,y},\re)\,dy\ =\
+\infty\hspace{2em}\mbox{for a.e.  }\xi\in\re^n,$$
hence both sides of (\ref{rel:fsec-f}) are equal to $+\infty$.
Indeed, if this is not the case, then we can find a set of linearly
independent vectors $\{\xi_1,\ldots,\xi_n\}$ such that (\ref{section4})
holds true for every $i\in\{1,\ldots,n\}$, hence $u\in \mathit{GSBV}(\re^n)$
(by statement (ii)), which is impossible.
\qed

\begin{rmk}
	\begin{em}
		If $\eta(\xi)$ is radial, \ie\ depends only on $|\xi|$, then
		the function $\left(S\varphi\right)(z)$ defined in
		(\ref{defn:hatphi}) depends only on $|z|$.  In particular, if
		$\varphi(z)=|z|^p$ (and $\eta$ is radial), then
		$\left(S\varphi\right)(z)= c_{p,n}\omega^{-1}|z|^p$, where
			\begin{equation}
				c_{p,n}:= \int_{S^{n-1}}\left|\langle
				v,e_{1}\rangle\right|^p\,d\hn(v),
				\hspace{2em}e_{1}:=(1,0,\ldots,0).
				\label{defn:cpn}
			\end{equation}
	\end{em}
\end{rmk}

Combining the results of \S\ \ref{sec:onedim} with equalities
(\ref{rel:fep-fepxi}), (\ref{rel:fepxi-fep-ndim}), and
(\ref{rel:fsec-f}) we can study the convergence of
$\{\mathcal{F}_{\ep}\}$.  For shortness' shake, we only give the
following result.

\begin{thm}
	\label{thm:ndim}
	Let $\varphi$ and $\psi$ be as in the lower semicontinuity Theorem
	\ref{thm:lsc}, and let $\eta\in L^{1}(\re^n)$ be a non-negative
	non-zero function.  Let $\{\phiep\}$ be a family of Borel
	functions satisfying assumptions (li1), (li2), and (Est) (cf.
	Theorem \ref{thm:liminf} and Theorem \ref{thm:estimate}) with
	$\varphi=\varphi_{\star}=\varphi^\star$ and
	$\psi=\psi_{\star}=\psi^\star$.  Finally, let
	$\omega:=\|\eta\|_{L^{1}(\re^n)}$, and let $S\varphi$ be the
	function defined in (\ref{defn:hatphi}).

	Then
	\begin{enumerate}
	\renewcommand{\labelenumi}{(\roman{enumi})\ }
		\item
		$\mathcal{F}_{\ep}(u,\re^n)\leq\omega\,
		\mathcal{F}_{S\varphi,\psi}(u,\re^n)$ for every $u\in
		L^{1}_{loc}(\re^n)$, and every $\ep>0$;

		\item $\{\mathcal{F}_{\ep}(u,\re^n)\}$ pointwise converges to
		$\omega\,\mathcal{F}_{S\varphi,\psi}(u,\re^n)$;

		\item $\omega\,\mathcal{F}_{S\varphi,\psi}(u,\re^n)$ is the
		$\Gamma^{\mbox{-}}$-limit of $\{\mathcal{F}_{\ep}(u,\re^n)\}$
		in $L^{1}_{loc}(\re^n)$.
	\end{enumerate}
\end{thm}
\prf
Let us prove statement (i).  By (\ref{rel:fep-fepxi}),
(\ref{rel:fepxi-fep-ndim}), and (\ref{rel:fsec-f}) we have that
\begin{eqnarray*}
	\mathcal{F}_{\ep}(u,\re^n) & = &
	\int_{\re^{n}}F_{\ep,\xi}(u,\re^n)\,\eta(\xi)\,d\xi \\
	\noalign{\vs}
	 & = & \int_{\re^n}^{}\left(\int_{\langle\xi\rangle^{\bot}}
	 F_{\ep|\xi|}(u_{\xi,y},\re)\,dy\right)\eta(\xi)\,d\xi \\
	\noalign{\vs}
	 & \leq & \int_{\re^n}^{}\left(\int_{\langle\xi\rangle^{\bot}}
	 \mathcal{F}_{\varphi,\psi}(u_{\xi,y},\re)\,dy\right)\eta(\xi)\,d\xi\
	 =\ \omega\,\mathcal{F}_{S\varphi,\psi}(u,\re^n).
\end{eqnarray*}

In order to complete the proof, it remains to show that
$$\liminf_{n\to\infty}\,\mathcal{F}_{\ep_{n}}(u_{n},\re^n)\ \geq\
\omega\,\mathcal{F}_{S\varphi,\psi}(u,\re^n)$$
for every sequence $\{\ep_{n}\}\to 0^{+}$, and every sequence
$\{u_{n}\}\to u$ in $L^{1}_{loc}(\re^n)$. Since
$$\{(u_{n})_{\xi,y}\}\lra u_{\xi,y}\hspace{2em} \mbox{in }L^1_{{\rm
loc}}(\re)$$
for a.e.  $y\in \langle\xi\rangle^{\bot}$, by (\ref{rel:fep-fepxi}),
(\ref{rel:fepxi-fep-ndim}), and Fatou's Lemma we have that:
\begin{eqnarray*}
	\liminf_{n\to\infty}\,\mathcal{F}_{\ep_{n}}(u_{n},\re^n) & = &
	\liminf_{n\to\infty}
	\int_{\re^{n}}F_{\ep_{n},\xi}(u_n,\re^n)\,\eta(\xi)\,d\xi \\
	\noalign{\vs}
	 & = &
	 \liminf_{n\to\infty}\int_{\re^n}\left(\int_{\langle\xi\rangle^{\bot}}
	 F_{\ep_{n}|\xi|}((u_n)_{\xi,y},\re)\,dy\right)\eta(\xi)\,d\xi \\
	\noalign{\vs}
	 & \geq &
	 \int_{\re^n}\left(\int_{\langle\xi\rangle^{\bot}}\liminf_{n\to\infty}
	 F_{\ep_{n}|\xi|}((u_n)_{\xi,y},\re)\,dy\right)\eta(\xi)\,d\xi \\
	\noalign{\vs}
	 & \geq & 	\int_{\re^n}\left(\int_{\langle\xi\rangle^{\bot}}
	 \mathcal{F}_{\varphi,\psi}(u_{\xi,y},\re)\,dy\right)\eta(\xi)\,d\xi \
	 =\ \omega\,\mathcal{F}_{S\varphi,\psi}(u,\re^n).
\end{eqnarray*}

This completes the proof.
\qed

\begin{rmk}\label{rmk:no-hp-phi}
	\begin{em}
		If in Theorem \ref{thm:ndim} we drop the assumption that
		$\varphi$ and $\psi$ satisfy (\ref{hp:phi}) and
		(\ref{hp:psi}), then the pointwise limit, the
		$\Gamma^{\mbox{-}}$-limit, and an upper estimate for
		$\mathcal{F}_\ep(u,\re^n)$ are given by the functional
		\begin{equation}
			\tilde{\mathcal{F}}(u,\re^n)\ :=\
			\int_{\re^n}\left(\int_{\langle\xi\rangle^{\bot}}
			\overline{\mathcal{F}_{\varphi,\psi}}
			(u_{\xi,y},\re)\,dy\right)\eta(\xi)\,d\xi,
			\label{ts:rmk}
		\end{equation}
		where $\overline{\mathcal{F}_{\varphi,\psi}}$ is the
relaxation of
		$\mathcal{F}_{\varphi,\psi}$ in the one-dimensional case.
	\end{em}
\end{rmk}

\begin{rmk}
	\begin{em}
		All the results of this section are true also in the
		particular case where $n=1$.  In this case
		$\langle\xi\rangle^{\bot}=\{0\}$ for every
		$\xi\in\re\setminus\{0\}$, and therefore many formulas
		containing integrations over $\langle\xi\rangle^{\bot}$ may be
		simplified.  Moreover (\ref{defn:hatphi}) reduces to
		$$\left(S\varphi\right)(z)\ =\ \frac{1}{\omega}
		\int_{\re}\varphi(|z|)\eta(\xi)\,d\xi\ =\ \varphi(|z|).$$
	\end{em}
\end{rmk}

\begin{rmk}
	\begin{em}
		All the results of this section (and in particular Theorem
		\ref{thm:ndim}) can be generalized to an arbitrary open set
		$\Omega\sse\re^n$.  In this case the natural generalization of
		(\ref{defn:fep-ndim-bis}) are the functionals
		\begin{equation}
			\mathcal{F}_{\ep}(u,\Omega)\ :=\
	\frac{1}{\ep^n}\int_{\mathrm{vis}(\Omega)}\varphi_{|y-x|}
			\left(\frac{|u(y)-u(x)|}{|y-x|}\right)
			\eta\left(\frac{y-x}{\ep}\right)\,dy\,dx,
		\label{defn:fep-ndim-omega}
		\end{equation}
		where
		$$\mathrm{vis}(\Omega):=\left\{(x,y)\in\Omega\times\Omega
		\tc\foralll t\in [0,1],\ tx+(1-t)y\in\Omega\right\}$$
		is the set of all pairs in $\Omega\times\Omega$ which ``see
		each other''.

		If $\Omega=\re^n$, then (\ref{defn:fep-ndim-bis}) can be
		written as (\ref{defn:fep-ndim-omega}) with the substitution
		$x+\ep\xi=y$.  The restriction of the integration to
		$\mathrm{vis}(\Omega)$, instead of $\Omega\times\Omega$, makes
		this construction to work on every open set $\Omega$,
		without any assumption on the regularity of the boundary
(see the
		discussion in \cite[section 7]{gobbino:ms}).
	\end{em}
\end{rmk}

\setcounter{equation}{0}
\section{Compactness}
\label{sec:compactness}
In this section we prove the following compactness result.
\begin{thm}\label{thm:cpt}
	Let $\{\phiep\}$ be a family of Borel functions such that
	\begin{description}
		 \item[\textnormal{(Cpt1)}] for every $M>0$ there exist
		 $H_{M}>0$ and $K_{M}\geq 0$ such that $$\phiep(r)\ \geq\
		 H_{M}\,r-K_{M}\hspace{2em}\foralll
		r\in[0,M/\ep];$$

		\item[\textnormal{(Cpt2)}]  $\phiep$ is nondecreasing for
every
		$\ep>0$, and
		$$\varphi_{(k+1)\ep}\left(\frac{A+B}{k+1}\right)\ \leq\
		\frac{1}{k+1}\,\phiep(A)+\frac{k}{k+1}\,
		\varphi_{k\ep}\left(\frac{B}{k}\right)$$
		for every $A\geq 0$, $B\geq 0$, $\ep>0$,
$k\in\n\setminus\{0\}$.
	\end{description}

	Let $\eta\in L^{1}(\re^n)$ be a non-negative function such that
	$\{\xi\in\re^n\tc\eta(\xi)>c\}$ has non-empty interior for some $c>0$.
	Let $\{\mathcal{F}_{\ep}\}$ be the functionals introduced in
	(\ref{defn:fep-ndim-bis}),
	and let $\{u_{\ep}\}\sse L^{\infty}(\re^{n})$ be such that
	\begin{equation}
		\sup_{\ep>0}\,\{{\cal
		F}_{\ep}(u_{\ep},\re^n)+\|u_{\ep}\|_{\infty}\}<+\infty.
		\label{hp:cpt}
	\end{equation}

	Then there exist $\{\ep_{k}\}\nlra 0^{+}$ and $u\in
\mathit{GSBV}(\re^{n})$
	such that
	$$\{u_{\ep_{k}}\}\lra u\hspace{2em}{\it in\ } L^{1}_{loc}(\re^{n}).$$
\end{thm}

\begin{rmk}
	\begin{em}
		If for some $p>1$ the inequality in (Cpt1) is replaced by
		$\phiep(r)\geq H_{M}r^p-K_{M}$, then the compactness result
holds
		true in $L^p_{loc}(\re^n)$.
	\end{em}
\end{rmk}

\subsection{Proofs}
In order to prove Theorem \ref{thm:cpt}, let us introduce some notations.
Let us assume that $\eta\in C^{1}(\re^n)$ is a non-zero non-negative
function with compact support, and let $R>0$ be such that
$$|\xi|\geq R\ \Longrightarrow\ \eta(\xi)=0.$$

Moreover we set
$$\omega_{0}:=\|\eta\|_{L^{1}(\re^n)},
\hspace{3em}\omega_{1}:=\|\nabla\eta\|_{L^{1}(\re^n)},$$
and we denote by $C^{\delta}u$ the convolutions
\begin{equation}
	C^{\delta}u(x)\ :=\ \frac{1}{\omega_{0}}\int_{\re^{n}}u(x+\delta\xi)\,
	\eta(\xi)\,d\xi,
	\label{defn:conv}
\end{equation}
defined for every $u\in L^{\infty}(\re^{n})$, and every $\delta>0$.

In a standard way it is possible to show that $C^{\delta}u\in
C^{1}(\re^{n})$ and moreover
	\begin{equation}
		\left\|C^{\delta}u\right\|_{\infty}\ \leq\
		\|u\|_{\infty},\hspace{3em} \left\|\nabla
		C^{\delta}u\right\|_{\infty}\ \leq\
		\frac{\omega_{1}}{\delta}\,\|u\|_{\infty}.
		\label{normcu}
	\end{equation}

We now need two technical lemmata.
\begin{lemma}\label{lemma:cpt1}
	Let $\{\phiep\}$ be a family of Borel functions satisfying (Cpt1).

	Then for all $u\in L^{\infty}(\re^{n})$, $\delta>0$,
	$A\subset\subset\re^{n}$, we have that
	$$\left\|C^{\delta}u-u\right\|_{L^{1}(A)}\ \leq\
	\frac{R}{H_{M}\omega_{0}}
	\left(\omega_{0}K_{M}|A|+\mathcal{F}_{\delta}(u,\re^n)\right)\delta$$
	for every $M\geq 2\|u\|_{\infty}$.
\end{lemma}
\prf
Applying (Cpt1) with
$$\ep=\delta|\xi|,\hspace{2em}
r=\frac{|u(x+\delta\xi)-u(x)|}{\delta|\xi|},$$
we have that
\begin{eqnarray*}
	|u(x+\delta\xi)-u(x)| & = &
	\delta|\xi|\frac{|u(x+\delta\xi)-u(x)|}{\delta|\xi|} \\
\noalign{\vs}
	 & \leq\ &
	 \frac{\delta|\xi|}{H_{M}}\left\{K_{M}+\varphi_{\delta|\xi|}\left(
	 \frac{|u(x+\delta\xi)-u(x)|}{\delta|\xi|}\right)\right\}
\end{eqnarray*}
for every $M\geq 2\|u\|_{\infty}$.  Therefore
\begin{eqnarray*}
	\lefteqn{\hspace{-2em}\left\|C^{\delta}u-u\right\|_{L^{1}(A)}\ \leq\
	\frac{1}{\omega_{0}}\int_{A\times\re^{n}}
	|u(x+\delta\xi)-u(x)|\,\eta(\xi)\,dx\,d\xi} \\
\noalign{\vs}
	 & \leq & \frac{\delta R}{H_{M}\omega_{0}} \int_{\re^n}
	 \int_{A} \left\{K_{M}+\varphi_{\delta|\xi|}\left(
	 \frac{|u(x+\delta\xi)-u(x)|}{\delta|\xi|}\right)\right\}dx\,
	 \eta(\xi)\,d\xi \\
\noalign{\vs}
	 & \leq &  \frac{\delta R}{H_{M}\omega_{0}} \int_{\re^n}\left(
	 K_{M}|A|+F_{\delta,\xi}(u,\re^n)\right)\,\eta(\xi)\,d\xi \\
\noalign{\vs}
	 & = & \frac{R}{H_{M}\omega_{0}}
	\left(\omega_{0}K_{M}|A|+\mathcal{F}_{\delta}(u,\re^n)\right)\delta.
\qedbis
\end{eqnarray*}

\begin{lemma}\label{lemma:cpt2}
	Let $\{\phiep\}$ be a family of Borel functions satisfying (Cpt2).

	Then for all $\delta>0$, $k\in\n\backslash\{0\}$, $u\in
	L^{\infty}(\re^{n})$ we have that
	\begin{equation}
		{\cal F}_{k\delta}(u,\re^n)\ \leq\ {\cal F}_{\delta}(u,\re^n).
		\label{ts:lemma-cpt3}
	\end{equation}
\end{lemma}
\prf
Let us argue by induction. If $k=1$, thesis is trivial. Let us
assume that (\ref{ts:lemma-cpt3}) holds true for some $k\geq 1$.
Applying (Cpt2) with $\ep=\delta|\xi|$, and
$$A=\frac{u(x+(k+1)\delta|\xi|)-u(x+k\delta|\xi|)}{\delta|\xi|}, \hspace{1em}
B=\frac{u(x+k\delta|\xi|)-u(x)}{\delta|\xi|},$$
it follows that
\begin{eqnarray*}
	\lefteqn{\hspace{-3em}\varphi_{(k+1)\delta|\xi|} \left(
	\frac{(u(x+(k+1)\delta|\xi|)-u(x))}{(k+1)\delta|\xi|}\right)\leq}\\
\noalign{\vs}
	 & \leq & \frac{1}{k+1}\,\varphi_{\delta|\xi|}\left(
\frac{(u(x+(k+1)\delta|\xi|)-u(x+k\delta|\xi|))}{\delta|\xi|}\right)+ \\
\noalign{\vs}
	 & & +\frac{k}{k+1}\,\varphi_{k\delta|\xi|}\left(
	 \frac{(u(x+k\delta|\xi|)-u(x))}{k\delta|\xi|}\right).
\end{eqnarray*}

Multiplying by $\eta(\xi)$ and integrating in $(x,\xi)$ over
$\re^{n}\times\re^{n}$, by the inductive hypothesis we obtain that
$${\cal F}_{(k+1)\delta}(u,\re^n)\ \leq\ \frac{1}{k+1}\,{\cal
F}_{\delta}(u,\re^n)\,+\, \frac{k}{k+1}\,{\cal F}_{k\delta}(u,\re^n)\
\leq\ {\cal F}_{\delta}(u,\re^n),$$
and this completes the induction.
\qed

\noindent\textsc{Proof of Theorem \ref{thm:cpt}}.

Up to replacing $\eta$ by a smaller function, we can assume that
$\eta$ belongs to $C^{1}(\re)$ and has compact support (this is the
point where we use our assumptions on $\eta$).
Now we argue as in the case of the Mumford-Shah functional.  We show
that $\{u_{\ep_{j}}\}$ is relatively compact in $L^{1}(A)$ for every
sequence $\{\ep_{j}\}\nlra 0^{+}$ and every $A\subset\subset\re^{n}$.
To this end we set for every $\sigma>0$
$$K_{\sigma}:=\left\{C^{\ep_{j}[\frac{\sigma}{\ep_{j}}]}
u_{\ep_{j}}\tc \ep_{j}\leq\frac{\sigma}{2}\right\}\cup
\bigcup_{\ep_{j}>\frac{\sigma}{2}}\{u_{\ep_{j}}\},$$
and we show that
\begin{enumerate}
	\renewcommand{\labelenumi}{(\roman{enumi})\ }
	\item  $K_{\sigma}$ is relatively compact in $L^{1}(A)$;

	\item  for all $j\in\n$, there exists $v_{j}\in K_{\sigma}$ such
	that $\|u_{\ep_{j}}-v_{j}\|_{L^{1}(A)}\leq N\,\sigma$, where $N$
	does not depend on $j$ and $\sigma$.
\end{enumerate}

This proves that the sequence $\{u_{\ep_{j}}\}$ is totally bounded,
hence relatively compact, in $L^{1}(A)$.

Let us show that $K_{\sigma}$ satisfies (i). Since there is only a
finite number of $\ep_{j}>\sigma/2$, it suffices to
show that
$$\tilde{K}_{\sigma}=\left\{C^{\ep_{j}[\frac{\sigma}{\ep_{j}}]}
u_{\ep_{j}}\tc \ep_{j}\leq\frac{\sigma}{2}\right\}$$
is relatively compact in $L^{1}(A)$.  To this end, let us remark that
$\tilde{K}_{\sigma}\sse C^{1}(A)$, and since
$\ep_{j}[\frac{\sigma}{\ep_{j}}]\geq\frac{\sigma}{2}$ by
(\ref{normcu}) we have that
$$\left\|C^{\ep_{j}[\frac{\sigma}{\ep_{j}}]}u_{\ep_{j}}\right\|_{\infty}\
\leq\ \|u_{\ep_{j}}\|_{\infty},\hspace{3em} \left\|\nabla
C^{\ep_{j}[\frac{\sigma}{\ep_{j}}]} u_{\ep_{j}}\right\|_{\infty}\
\leq\ \frac{2\omega_{1}}{\sigma}\,\|u_{\ep_{j}}\|_{\infty}.$$

By Ascoli's Theorem, $\tilde{K}_{\sigma}$ is relatively
compact in $C^{0}(A)$, hence in $L^{1}(A)$.

Let us show that $K_{\sigma}$ satisfies (ii) with
$$N:=\frac{R}{H_{M}\omega_{0}}\left(\omega_{0}K_{M}|A|+
\sup_{\ep>0}\mathcal{F}_{\ep}(u_{\ep},\re^n)\right),\hspace{2em}
M:=2\sup_{\ep>0}\|u_{\ep}\|_{\infty}.$$

If $\ep_{j}>\sigma/2$ we can simply take $v_{j}=u_{\ep_{j}}$. If
$\ep_{j}\leq\sigma/2$ we can take
$v_{j}=C^{\ep_{j}[\frac{\sigma}{\ep_{j}}]}u_{\ep_{j}}$. Indeed, by
Lemma\,\ref{lemma:cpt1} and Lemma\,\ref{lemma:cpt2}, we have that
\begin{eqnarray*}
	\left\|C^{\ep_{j}[\frac{\sigma}{\ep_{j}}]}
	u_{\ep_{j}}-u_{\ep_{j}}\right\|_{L^{1}(A)} & \leq &
	\frac{R}{H_{M}\omega_{0}} \left(\omega_{0}K_{M}|A|+
	\mathcal{F}_{\ep_{j}[\frac{\sigma}{\ep_{j}}]}
	(u_{\ep_{j}},\re^n))\right)\sigma \\
\noalign{\vs}
	& \leq & \frac{R}{H_{M}\omega_{0}}
	\left(\omega_{0}K_{M}|A|+
	\mathcal{F}_{\ep_{j}}(u_{\ep_{j}},\re^n))\right)\sigma\ \leq\
	N\sigma.
\end{eqnarray*}

By (\ref{hp:cpt}) and the liminf inequality in the definition of
$\Gamma^{\mbox{-}}$-convergence, any limit point of $\{u_{\ep}\}$
satisfies $\mathcal{F}_{\varphi,\psi}(u,\re^n)<+\infty$, hence
belongs to $GSBV(\re^n)$.
\qed

\setcounter{equation}{0}
\section{Approximation of Free Discontinuity Problems}
\label{sec:approximation}
In this section we prove that large classes of functionals like
(\ref{defn:f-phi-psi}) can be approximated by non-local functionals
of the form (\ref{defn:fep-ndim-bis}). To this end, we need
the following definition.

\begin{defn}\label{defn:sectionable}
	\begin{em}
		We say that an increasing convex function
		$\varphi:[0,+\infty[\to[0,+\infty]$ is \emph{sectionable} in
		$\,\re^n$ if there exists a convex function
		$\overline{\varphi}:[0,+\infty[\to[0,+\infty]$ such that
		\begin{equation}
			\varphi(|\alpha|)=\frac{1}{\hn(S^{n-1})}\int_{S^{n-1}}
	\overline{\varphi}\left(|\langle\alpha,v\rangle|\right) d\hn(v)
			\hspace{2em}\foralll\alpha\in\re^n.
			\label{defn:phi-sect}
		\end{equation}
	\end{em}
\end{defn}

\begin{rmk}
	\begin{em}
		The following properties of sectionable functions are an
		immediate consequence of the above definition.
	\begin{itemize}
		\item Every convex function
		$\varphi:[0,+\infty[\to[0,+\infty]$ is sectionable in $\re$,
		and $\overline{\varphi}=\varphi$.

		\item For every real number $p\geq 1$, the function
		$\varphi(r)=r^p$ is sectionable in $\re^n$ for every $n$.
		Indeed, (\ref{defn:phi-sect}) is satisfied with
		$\overline{\varphi}(r)=(c_{p,n})^{-1}r^p$, where $c_{p,n}$
		is the constant introduced in (\ref{defn:cpn}).

		\item The class of sectionable functions is additively
closed.
		Moreover, if $\varphi$ is sectionable, and $\lambda>0$ is a
		constant, then $\lambda\varphi$ is sectionable.

		\item The class of sectionable functions is closed by monotone
		convergence in the following sense: if $\{\varphi_{n}\}$ is a
		sequence of sectionable functions, and
		$\varphi_{n+1}(r)\geq\varphi_{n}(r)$ for every $n\in\n$ and
		every $r\geq 0$, then $\sup\varphi_{n}$ is sectionable.  In
		this way we can show, for example, that $\varphi(r)=e^{r^{2}}$
		is sectionable.

		\item Every sectionable function $\varphi$ is the supremum of
		an increasing sequence of sectionable \emph{finite} functions
		(it is enough to approximate $\overline{\varphi}$ with an
		increasing sequence of finite convex functions).

		\item If $\varphi$ is sectionable in $\re^n$ and satisfies
		(\ref{hp:phi}), then also $\overline{\varphi}$ satisfies
		(\ref{hp:phi}).

		\item It can be proved that $\varphi(r):=\max\{0,r-1\}$ is
		\emph{not} sectionable in $\re^n$ for every $n>1$.
		\end{itemize}
	\end{em}
\end{rmk}

The following is the main result of this paper.
\begin{thm}
	\label{thm:main}
	Let $\varphi$ and $\psi$ be as in the lower semicontinuity Theorem
	\ref{thm:lsc}, and let $\eta\in L^{1}(\re^n)$ be a non-negative
	radial function such that $\{\xi\in\re^n\tc\eta(\xi)>c\}$ has
	non-empty interior for some $c>0$.  Let us assume that $\varphi$
	is sectionable in $\re^n$.

	Then there exists a family $\{\phiep\}$ such that, defining
	$\left\{\mathcal{F}_{\ep}\right\}$ as in (\ref{defn:fep-ndim-bis}), we
	have that
		\begin{enumerate}
		\renewcommand{\labelenumi}{(C\arabic{enumi})\ }
		\item

	$\mathcal{F}_{\ep}(u,\re^n)\leq\mathcal{F}_{\varphi,\psi}(u,\re^n)$
		for all $u\in L^{1}_{loc}(\re^n)$, and all $\ep>0$;

		\item $\{\mathcal{F}_{\ep}(u,\re^n)\}$ converges to
		$\mathcal{F}_{\varphi,\psi}(u,\re^n)$ for all $u\in
		L^{1}_{loc}(\re^n)$;

		\item  $\mathcal{F}_{\varphi,\psi}(u,\re^n)$ is the
		$\Gamma^{\mbox{-}}$-limit of
$\{\mathcal{F}_{\ep}(u,\re^n)\}$ in
		$L^{1}_{loc}(\re^n)$;

		\item  	if $\{u_{\ep}\}\sse L^{\infty}(\re^{n})$ and
		$$\sup_{\ep>0}\left\{\mathcal{F}_{\ep}(u_{\ep},\re^n)
		+\|u_{\ep}\|_{\infty}\right\}\ <\ +\infty,$$
	then there exist $\{\ep_{k}\}\nlra 0^{+}$ and $u\in
	\mathit{GSBV}(\re^{n})$ such that
	$$\left\{u_{\ep_{k}}\right\}\lra
	u\hspace{2em}{\it in\ } L^{1}_{loc}(\re^{n}).$$
	\end{enumerate}
\end{thm}

\begin{rmk}
         \begin{em}
		For simplicity's sake we developed our theory under the
		assumptions of Theorem \ref{thm:lsc}, as stated in \cite{
		ambrsc, ambrarch}. However, it is well known that Theorem
\ref{thm:lsc}
		holds true also when the assumption
		``$\psi$ is concave'' is relaxed to ``$\psi$ is
		sub-additive and lower semi-continuous'' (see \eg\
\cite[Theorem
		2.10]{braides:survey}). In the same way, throughout all
this paper
		(hence in Theorem \ref{thm:main} above), we can modify the
		concavity assumptions on $\psi$ to sub-additivity and lower
		semi-continuity (but some proofs may become longer!).

         \end{em}
\end{rmk}

\begin{rmk}
	\begin{em}
		The family $\{\phiep\}$ given by Theorem \ref{thm:main} is
		clearly not unique.  In our proof, $\phiep$ will be defined as
		the minimum of a family of functions.  This construction is
		convenient from the theoretic point of view, but often it is
		difficult to give an explicit expression of this minimum.  For
		this reason, in many applications it may be useful to find
		other families with a simpler analytic expression, and then
		prove the convergence case-by-case using Theorem
		\ref{thm:ndim} and Theorem \ref{thm:cpt} (cf.  the examples in
		\S\ \ref{sec:examples}).
	\end{em}
\end{rmk}

\subsection{Proofs}
In this subsection we prove Theorem \ref{thm:main}. To this end, we
need two lemmata about real functions.

\begin{lemma}\label{lemma:f-conv-g}
	Let $f:[0,+\infty[\to [0,+\infty[$ be a non-decreasing convex
	function, and let $g:[0,+\infty]\to [0,+\infty[$ be a
	non-decreasing concave function such that $g(0)=0$.  Let us set
	\begin{equation}
		\mu(r):=\min\{f(l)+g(r-l)\tc l\in[0,r]\}.
		\label{defn:mu}
	\end{equation}

	Then
	\begin{enumerate}
	\renewcommand{\labelenumi}{(\roman{enumi})\ }
		\item  $\mu$ is continuous and non-decreasing;

		\item  there exists $\overline{r}\geq 0$ such that $\mu$ is
convex
		in $[0,\overline{r}]$ and concave in $[\overline{r},+\infty[$.
	\end{enumerate}
\end{lemma}
\prf
Let us first remark that our assumptions imply the continuity of $f$
in $[0,+\infty[$, and the continuity of $g$ in $]0,+\infty[$ (but not
the continuity of $g$ in $[0,+\infty[\,$!). In any case, $g$ is at
least lower semicontinuous in $[0,+\infty[$, and therefore the
``$\min$'' in (\ref{defn:mu}) is attained. Moreover, by definition of
$\mu$ we have that
\begin{equation}
	\mu(r)\leq f(r),\hspace{2em}\foralll r\geq 0.
	\label{estim:mu-leq-f}
\end{equation}

We claim that $\mu$ satisfies (ii) with
$$\overline{r}:=\sup\{r\geq 0\tc\mu(r)=f(r)\}.$$

\noindent\emph{Step 1.} We prove that $\mu$ is non-decreasing. To this
end, let $s>r$, and let $l\in[0,s]$ be such that $\mu(s)=f(l)+g(s-l)$.

If $l\in[0,r]$, by the monotonicity of $g$ we have that
$$\mu(s)\ =\ f(l)+g(s-l)\;\geq \; f(l)+g(r-l)\;\geq\;\mu(r).$$

If $l\in]r,s]$, by the monotonicity of $f$ and (\ref{estim:mu-leq-f})
it follows that
$$\mu(s)\ =\ f(l)+g(s-l)\;\geq \;f(r)\;\geq\;\mu(r).$$

In any case, we have proved that $\mu(s)\geq\mu(r)$.\vs

\noindent\emph{Step 2.} We show that
\begin{equation}
	\mu(r)=f(r),\hspace{2em}\foralll r\in[0,\overline{r}],
	\label{ts:mu=f}
\end{equation}
and therefore $\mu$ is convex in  $[0,\overline{r}]$.

Indeed, let us assume by contradiction that $\mu(r_{*})<f(r_{*})$ for
some $r_{*}\in[0,\overline{r}[$. Then, there exists $l<r_{*}$ such
that
$$f(l)+g(r_{*}-l)<f(r_{*}).$$

Now let us consider the function $\gamma:[l,+\infty[\to\re$ defined by
$$\gamma(t):=f(t)-g(t-l)-f(l).$$

Since $\gamma$ is convex, $\gamma(l)=0$, and $\gamma(r_{*})>0$, then
necessarily $\gamma(t)>0$ for every $t\geq r_{*}$.  Therefore
$$\mu(r)\ \leq\ f(l)+g(r-l)\ <\ f(r)\hspace{2em}\foralll r\geq r_{*},$$
which contradicts the definition of $\overline{r}$. This proves that
$\mu(r)=f(r)$ for every $r\in[0,\overline{r}[$. By the monotonicity
of $\mu$ and (\ref{estim:mu-leq-f}) it follows that
$$f(r)\ =\ \mu(r)\ \leq\ \mu(\overline{r})\ \leq f(\overline{r})$$
for every $r<\overline{r}$. Passing to the limit as
$r\to\overline{r}^-$, the proof of (\ref{ts:mu=f}) is complete.\vs

\noindent\emph{Step 3.} We prove that there exists
$\overline{l}\in[0,\overline{r}]$ such that
\begin{equation}
	f\left(\overline{l}\right)+g\left(r-\overline{l}\right)
	\ \leq\ f(r),\hspace{2em}\foralll
	r\geq\overline{r}.
	\label{ts:over-l}
\end{equation}

Indeed, let $\{r_{n}\}\to\overline{r}^{+}$ be any sequence, and for
each $n\in\n$, let $l_{n}\in[0,r_{n}]$ be such that
$$\mu(r_{n})\ =\ f(l_{n})+g(r_{n}-l_{n})\ <\ f(r_{n}),$$
where the inequality follows from the definition of
$\overline{r}$.

Up to subsequences, we can assume that
$\{l_{n}\}\to\overline{l}\in[0,\overline{r}]$. In order to prove that
(\ref{ts:over-l}) holds true, let us fix $r>\overline{r}$, and let us
consider the functions $\gamma_{n}:[l_{n},+\infty[\to\re$ defined by
$$\gamma_{n}(t)\ :=\ f(t)-g(t-l_{n})-f(l_{n}).$$

Since $\gamma_{n}$ is a convex function such that
$\gamma_{n}(l_{n})=0$ and $\gamma_{n}(r_{n})>0$, then necessarily
$\gamma_{n}(t)>0$ for every $t\geq r_{n}$. Since $r\geq r_{n}$ for $n$
large enough, it follows that $\gamma_{n}(r)>0$ for $n$ large enough.
Passing to the limit as $n\to+\infty$ we obtain that
$$f(r)-g(r-\overline{l})-f(\overline{l})\geq 0,$$
which is equivalent to (\ref{ts:over-l}).\vs

\noindent\emph{Step 4.} We prove that
\begin{equation}
	\mu(r)=\min\{f(l)+g(r-l)\tc l\in[0,\overline{r}]\},
	\hspace{2em}\foralll r\geq \overline{r}.
	\label{ts:restr-min}
\end{equation}

Indeed, if $r\geq l\geq\overline{r}$, then, using (\ref{ts:over-l})
with $r=l$, and the subadditivity of $g$, it follows that
$$f(l)+g(r-l)\ \geq\ f(\overline{l})+g(l-\overline{l})+g(r-l)\ \geq\
f(\overline{l})+g(r-\overline{l}).$$

This proves that for $r\geq \overline{r}$, in the minimum problem
(\ref{defn:mu}) we can consider only the values
$l\in[0,\overline{r}]$.\vs

\noindent\emph{Step 5.} By (\ref{ts:restr-min}) we have that for $r\geq
\overline{r}$, the function $\mu$ is the minimum of a \emph{fixed} family of
concave functions. This proves that $\mu$ is concave in
$[\overline{r},+\infty[$.\vs

\noindent\emph{Step 6.} In order to complete the proof of the lemma,
it remains to show that $\mu$ is continuous.

By (\ref{ts:mu=f}) the restriction of $\mu$ to $[0,\overline{r}]$ is
continuous.  Moreover, $\mu$ is continuous on $]\overline{r},+\infty[$
since it is concave in this region.  Therefore it remains to prove that
\begin{equation}
	\lim_{r\to\overline{r}^{+}}\mu(r)=\mu(\overline{r}).
	\label{ts:mu-cont}
\end{equation}

By the monotonicity of $\mu$, and (\ref{ts:over-l}), it follows that
$$f(\overline{r})\ =\ \mu(\overline{r})\ \leq\ \mu(r)\ \leq\
f(\overline{l})+g(r-\overline{l})\ \leq\ f(r)\hspace{3em}
\foralll r\geq\overline{r}.$$

Passing to the limit as $r\to\overline{r}^{+}$, (\ref{ts:mu-cont}) is
proved.\vs
\qed

\begin{lemma}\label{lemma:phiep}
	Let $\varphi:[0,+\infty[\to [0,+\infty[$ be a non-decreasing
	convex function not identically zero, and let $\psi:[0,+\infty]\to
	[0,+\infty[$ be a non-decreasing concave function such that
$\psi(0)=0$.  Let us set
	\begin{equation}
		\phiep(r)\ :=\ \min\left\{\varphi(l)+\frac{1}{\ep}\,
		\psi(\ep(r-l))\tc l\in[0,r]\right\}.
		\label{defn:phiep-phi-psi}
	\end{equation}
	for every $\ep>0$.
	Then
\begin{enumerate}
	\renewcommand{\labelenumi}{(\roman{enumi})\ }
	\item the family $\{\phiep\}$ satisfies (li1), (li2), (Est),
	(Cpt1) and (Cpt2).
\end{enumerate}

If moreover $\varphi$ and $\psi$ satisfy (\ref{hp:phi}) and
(\ref{hp:psi}), then
\begin{enumerate}
\setcounter{enumi}{1}
	\renewcommand{\labelenumi}{(\roman{enumi})\ }
	\item $\{\phiep(r)\}\to\varphi(r)$ uniformly on compact subsets of
	$[0,+\infty[$;

	\item $\{\ep\phiep(r/\ep)\}\to\psi(r)$ uniformly on compact subsets
	of $]0,+\infty[$.
\end{enumerate}

In particular
$$\varphi(r)=\gmliminf{\ep\to 0^{+}}\,\phiep(r),
\hspace{2em}\psi(r)=\gmliminf{\ep\to
0^{+}}\,\ep\phiep\left(\frac{r}{\ep}\right).$$
\end{lemma}
\prf\\
\noindent \emph{Proof of (i).} Properties (li1) and (li2) follow from
Lemma \ref{lemma:f-conv-g} applied with $f(r)=\varphi(r)$ and
$g(r)=\psi(\ep r)/\ep$. Property (Est) is a trivial consequence of
the definition (\ref{defn:phiep-phi-psi}).

Since $\varphi$ is convex and non-zero, then there exists $c>0$ and
$d\geq 0$ such that
$$\varphi(r)\ \geq\ cr-d\hspace{2em}\foralll r\geq 0.$$

Moreover, since $\psi$ is concave, then for every $M>0$ we have that
$$\psi(r)\ \geq\ \frac{\psi(M)}{M}\,r\hspace{2em}\foralll r\in[0,M].$$

We claim that $\phiep$ satisfies (Cpt1) with
$$H_{M}:=\min\left\{\frac{\psi(M)}{M},c\right\},\hspace{2em}K_{M}=d.$$

Indeed, for every $0\leq l\leq r\leq M/\ep$ we have that
$\ep(r-l)\leq M$, hence
$$\varphi(l)+\frac{1}{\ep}\,\psi(\ep(r-l))\ \geq\
H_{M}l-K_{M}+\frac{1}{\ep}\,H_{M}\ep(r-l)\ =\ H_{M}\,r-K_{M}.$$

This is equivalent to (Cpt1).

Now let us prove that $\{\phiep\}$ satisfies (Cpt2). Let $l_{A}\in
[0,A]$ and $l_{B}\in[0,B/k]$ be such that
$$\varphi_{\delta}(A)\ =\
\varphi(l_{A})+\frac{1}{\delta}\,\psi(\delta(A-l_{A})),$$
$$\varphi_{k\delta}\left(\frac{B}{k}\right)\ =\
\varphi(l_{B})+\frac{1}{k\delta}\,\psi
\left(k\delta\left(\frac{B}{k}-l_{B}\right)\right).$$

Then
$$l_{C}\ :=\ \frac{l_{A}+kl_{B}}{k+1}\in
\left[0,\frac{A+B}{k+1}\right],$$
hence, by the convexity of $\varphi$ and the subadditivity of $\psi$,
it follows that
\begin{eqnarray*}
	\lefteqn{\hspace{-2em}\varphi_{(k+1)\delta}\left(\frac{A+B}{k+1}
\right) \ \leq\
	\varphi\left(l_{C}\right)
	+\frac{1}{(k+1)\delta}\,\psi
	\left((k+1)\delta\left(\frac{A+B}{k+1}-
	l_{C}\right)\right)}  \\
	\noalign{\vs}
	 & = & \varphi\left(\frac{l_{A}}{k+1}+\frac{k}{k+1}l_{B}\right)
	 +\frac{1}{(k+1)\delta}\,\psi\left(\delta(A-l_{A})
	 +k\delta\left(\frac{B}{k}-l_{B}\right)\right)  \\
	\noalign{\vs}
	 & \leq &
\frac{1}{k+1}\,\varphi(l_{A})+\frac{k}{k+1}\,\varphi(l_{B}) \\
	\noalign{\vs}
	 & & +\frac{1}{(k+1)\delta}\,\psi(\delta(A-l_{A}))
	 +\frac{k}{(k+1)}\,\frac{1}{k\delta}\,\psi
	 \left(k\delta\left(\frac{B}{k}-l_{B}\right)\right) \\
	\noalign{\vs}
	 & = & \frac{1}{k+1}\,\varphi_{\delta}(A)+ \frac{k}{k+1}
	 \,\varphi_{k\delta}\left(\frac{B}{k}\right).
\end{eqnarray*}

\noindent \emph{Proof of (ii).} Let $r\geq 0$.  Setting $l=r$ in
(\ref{defn:phiep-phi-psi}) we have that
\begin{equation}
	\phiep(r)\leq\varphi(r).
	\label{estim:phiep1}
\end{equation}

Moreover, let $l_{\ep}\in[0,r]$ be such that
\begin{equation}
	\phiep(r)=\varphi(l_{\ep})+\frac{1}{\ep}\,\psi(\ep(r-l_{\ep})).
	\label{defnlep1}
\end{equation}

We claim that $\{l_{\ep}\}\to r$. Indeed, let us assume by contradiction
that there exists a sequence $\{\ep_{n}\}\to 0^{+}$ such that
$\left\{r-l_{\ep_{n}}\right\}\to\alpha>0$. Then, since
$\left\{\ep_{n}(r-l_{\ep_{n}})\right\}\to 0$, by (\ref{defnlep1})
and (\ref{hp:psi}) it follows that
\begin{eqnarray*}
	\lim_{n\to +\infty}\,\varphi_{\ep_{n}}(r) & = & \lim_{n\to +\infty}
	\,\varphi(l_{\ep_{n}})+\frac{1}{\ep_{n}}\,
	\psi(\ep_{n}(r-l_{\ep_{n}})) \\
	\noalign{\vs}
	 & \geq & \lim_{n\to
	+\infty}\,(r-l_{\ep_{n}})\,\frac{\psi(\ep_{n}(r-l_{\ep_{n}}))}
	{\ep_{n}(r-l_{\ep_{n}})}\ =\ +\infty,
\end{eqnarray*}
which is impossible because of (\ref{estim:phiep1}).  Since
$\{l_{\ep}\}\to r$, then by (\ref{defnlep1}) it follows that
$$\liminf_{\ep\to 0^{+}}\,\phiep(r)\ \geq\ \liminf_{\ep\to
0^{+}}\,\varphi(l_{\ep})\ =\ \varphi(r),$$
which, together with (\ref{estim:phiep1}), proves that
$\{\phiep(r)\}\to\varphi(r)$ for all $r\geq 0$.
Since $\phiep$ and $\varphi$ are continuous increasing functions,
uniform convergence on compact subsets follows from pointwise
convergence.\vs

\noindent\emph{Proof of (iii).} Let $r> 0$.  Setting $l=0$ in
(\ref{defn:phiep-phi-psi}) we have that
\begin{equation}
	\ep\phiep\left(\frac{r}{\ep}\right)\leq\psi(r)+\ep\varphi(0).
	\label{estim:phiep2}
\end{equation}

Moreover, let $l_{\ep}\in[0,r/\ep]$ be such that
\begin{equation}
	\phiep\left(\frac{r}{\ep}\right)=
	\varphi(l_{\ep})+\frac{1}{\ep}\,\psi(r-\ep l_{\ep}).
	\label{defnlep2}
\end{equation}

We claim that $\{\ep l_{\ep}\}\to 0$.  Indeed, let us assume by
contradiction that there exists a sequence $\{\ep_{n}\}\to 0^{+}$ such
that $\left\{\ep_{n}l_{\ep_{n}}\right\}\to\alpha\in]0,r]$.  Then
$\{l_{\ep_{n}}\}\to +\infty$ hence, by (\ref{defnlep2}) and
(\ref{hp:phi}),
\begin{eqnarray*}
	\lim_{n\to
	+\infty}\,\ep_{n}\varphi_{\ep_{n}}\left(\frac{r}{\ep_{n}}\right) &
	= & \lim_{n\to +\infty}\,
	\ep_{n}\varphi(l_{\ep_{n}})+\psi(r-\ep_{n}l_{\ep_{n}}) \\
	 & \geq & \lim_{n\to +\infty}\ep_{n}l_{\ep_{n}}\,
\frac{\varphi(l_{\ep_{n}})}{l_{\ep_{n}}}=+\infty,
\end{eqnarray*}
which is impossible because of (\ref{estim:phiep2}). Since $\{\ep l_{\ep}\}\to
0$, then by (\ref{defnlep2}) it follows that
$$\liminf_{\ep\to 0^{+}}\,\ep\phiep\left(\frac{r}{\ep}\right)\ \geq\
\liminf_{\ep\to 0^{+}}\,\psi(r-\ep l_{\ep})\ =\ \psi(r),$$
which, together with (\ref{estim:phiep2}), proves that
$\{\ep\phiep(r/\ep)\}\to\psi(r)$ for every $r>0$.  Since
$\ep\phiep(r/\ep)$ and $\psi(r)$ are continuous increasing functions,
uniform convergence on compact subsets follows from pointwise
convergence.
\qed

\noindent\textsc{Proof of Theorem \ref{thm:main}.}

Let $\eta(\xi)=J(|\xi|)$ for some $J:[0,+\infty[\to[0,+\infty[$.  Let
us assume that $\varphi$ and $\psi$ are finite.  In this case, we
extend $\psi$ to $[0,+\infty]$ by setting $\psi(0)=0$, and then, for
every $\ep>0$, we define
$$\phiep(r):=\frac{1}{\omega}\min
	\left\{\overline{\varphi}(l)+\frac{1}{\ep}\,\psi(\ep (r-l))\tc
	l\in[0,r]\right\},\hspace{2em}\foralll r\geq 0,$$
where $\overline{\varphi}$ is given by (\ref{defn:phi-sect}), and
$\omega:=\|\eta\|_{L^{1}(\re^n)}$.

By Lemma \ref{lemma:phiep} it
follows that $\{\phiep\}$ satisfies (li1), (li2) and (Est) with
$$\varphi_{\star}(r)=\varphi^{\star}(r)=
\frac{\overline{\varphi}(r)}{\omega},\hspace{2em}
\psi_{\star}(r)=\psi^{\star}(r)= \frac{\psi(r)}{\omega}.$$

Since using spherical coordinates in (\ref{defn:hatphi}) we have that
\begin{eqnarray*}
	\left(S(\overline{\varphi}/\omega)\right)(\alpha) & = &
	\frac{1}{\omega}\int_{\re^n}\frac{1}{\omega}\,\overline{\varphi}
	\left(|\langle\alpha,\xi/|\xi|\rangle|\right)\,\eta(\xi)\,d\xi  \\
	 \noalign{\vs}
	 & = &
	 \frac{1}{\omega^{2}}\int_{0}^{+\infty}\rho^{n-1}\,J(\rho)\,d\rho
	 \int_{S^{n-1}}\overline{\varphi}
	 (|\langle\alpha,v\rangle|)\,d\hn(v) \\
	 \noalign{\vs}
	 & = & \frac{1}{\omega^{2}}\int_{0}^{+\infty}\rho^{n-1}\,J(\rho)
	 \,\hn(S^{n-1})\,\varphi(|\alpha|)\,d\rho\ =\
	 \frac{1}{\omega}\,\varphi(|\alpha|),
\end{eqnarray*}
then statements (C1), (C2), and (C3) follow from Theorem \ref{thm:ndim}.

Moreover, by Lemma \ref{lemma:phiep} we have that $\{\phiep\}$
satisfies also (Cpt1) and (Cpt2), hence statement (C4) follows from
Theorem \ref{thm:cpt}.

If $\varphi$ and $\psi$ are not finite, then we first
approximate $\mathcal{F}_{\varphi,\psi}$ from below by functionals
$\{\mathcal{F}_{\varphi_{n},\psi_{n}}\}$, where $\varphi_{n}$ and
$\psi_{n}$ are finite functions satisfying the assumptions of this
theorem.  Arguing as before, we approximate the functionals
$\mathcal{F}_{\varphi_{n},\psi_{n}}$, and then we conclude the proof
by a diagonal argument.
\qed

\setcounter{equation}{0}
\section{Examples}
\label{sec:examples}
In this section we give some applications of the results proved in the
previous sections.  We apply Theorem \ref{thm:ndim} and Theorem
\ref{thm:cpt} in order to prove the convergence results (C1) through
(C4) of \S\ \ref{sec:introduction} for some special choices of
$\{\phiep\}$.

From now on, we assume that $J:[0,+\infty[\to[0,+\infty[$ is a
continuous function not identically zero such that
$$j_{\alpha}:=\int_{0}^{+\infty}\rho^{\alpha-1}J(\rho)\,d\rho$$
is finite for each real number $\alpha\geq 1$.  We also consider the
constants $c_{p,n}$ defined in (\ref{defn:cpn}).  In particular:
$c_{0,n}=\hn(S^{n-1})$, $J(|\xi|)\in L^{1}(\re^n)$, and
$\|\eta\|_{L^{1}(\re^n)}=c_{0,n}\,j_{n}$.

\begin{example} \label{ex:rp}
\begin{em}
	Let us consider the functionals
	$$\mathcal{F}_{\ep}(u)\ :=\ \frac{1}{\ep^{p}} \int_{\re^n\times\re^n}
	\left|u(x+\ep\xi)-u(x)\right|^{p}J(|\xi|)\,d\xi\,dx,$$
	with $p>1$.  Then $\left\{\mathcal{F}_{\ep}\right\}$ satisfies
	(C1) through (C4) of \S\ \ref{sec:introduction} with
	$$\mathcal{F}(u)\ :=\ \cases{{\displaystyle \lambda\int_{\re^n}
	|\nabla
	u(x)|^{p} d x} & \hspace{1em}if $u \in W^{1,p}_{loc}(\re^n)$, \cr
	\noalign{\vs}
	+\infty & \hspace{1em}if $u \in L^{1}_{loc}(\re^n)
	\setminus W^{1,p}_{loc}(\re^n),$ \cr}$$
	where $\lambda=c_{p,n}\,j_{p+n}$.

	Indeed the family $\left\{\mathcal{F}_{\ep}\right\}$ is a particular
	instance of (\ref{defn:fep-ndim-bis}) with
	$$\phiep(r)\ :=\ |r|^p,\hspace{4em} \eta(\xi)\ :=\ |\xi|^pJ(|\xi|).$$

	Since $\phiep$ satisfies (li1), (li2), (Est), (Cpt1), (Cpt2) with
	$$\varphi_{\star}(r)=\varphi^{\star}(r)=|r|^p,\hspace{4em}
	\psi_{\star}(r)=\psi^{\star}(r)=+\infty,$$
	then the results follow from Theorem \ref{thm:ndim} and Theorem
	\ref{thm:cpt}.
\end{em}
\end{example}

\begin{example} \label{ex:r1/p}
\begin{em}
	Let us consider the functionals
	$$\mathcal{F}_{\ep}(u)\ :=\
	\frac{1}{\ep} \int_{\re^n\times\re^n}
	\left|u(x+\ep\xi)-u(x)\right|^{1/p}J(|\xi|)\,d\xi\,dx,$$
	with $p>1$.  Then $\left\{\mathcal{F}_{\ep}\right\}$ satisfies
	(C1) through (C4) of \S\ \ref{sec:introduction} with
	$$\mathcal{F}(u)\ :=\ \cases{\displaystyle {\lambda\int_{S_{u}}
|u^{+} -
	u^{-}|^{1/p}\,d\hn} & \hspace{1em}if $u \in \mathit{GSBV}(\re^n)\cap
	L^{1}_{loc}(\re^n)$ and\cr
	\noalign{\vspace{-1mm}}
	 & \hspace{1em}$\nabla u(x) = 0$ for
	a.e.  $x\in\re^n$, \cr
	\noalign{\vs}
	+\infty & \hspace{1em}if $u
	\in L^{1}_{loc}(\re^n)\setminus \mathit{GSBV}(\re^n)$,\cr}$$
	where $\lambda=c_{0,n}\,j_{n+1/p}$.

	Indeed the family $\left\{\mathcal{F}_{\ep}\right\}$ is a
	particular instance of (\ref{defn:fep-ndim-bis}) with
	$$\phiep(r)\
	:=\ \frac{1}{\ep}\,|\ep r|^{1/p},\hspace{4em} \eta(\xi)\ :=\
	|\xi|^{1/p}J(|\xi|).$$

	Since $\phiep$ satisfies (li1), (li2), (Est), (Cpt1), (Cpt2) with
	$$\varphi_{\star}(r)=\varphi^{\star}(r)=+\infty,\hspace{4em}
	\psi_{\star}(r)=\psi^{\star}(r)=|r|^{1/p},$$
	then the results follow from Theorem \ref{thm:ndim} and Theorem
	\ref{thm:cpt}.
\end{em}
\end{example}

\begin{example} \label{ex:r}
\begin{em}
	Let us consider the functionals
	$$\mathcal{F}_{\ep}(u)\ :=\ \frac{1}{\ep} \int_{\re^n\times\re^n}
	\left|u(x+\ep\xi)-u(x)\right|J(|\xi|)\,d\xi\,dx.$$

	This is the limit case $p=1$ both of Example \ref{ex:rp} and of
	Example \ref{ex:r1/p}.  In this case (C1) through (C4) of \S\
	\ref{sec:introduction} are satisfied with
	$$\mathcal{F}(u)\ :=\ \lambda\left|Du\right|(\re^n),$$
	where $\lambda=c_{1,n}\,j_{n+1}$, and $\left|Du\right|(\re^n)$
	is defined as in (\ref{defn:du}).

	Indeed the family $\left\{\mathcal{F}_{\ep}\right\}$ is a particular
	instance of (\ref{defn:fep-ndim-bis}) with
	$$\phiep(r)\ :=\ |r|,\hspace{4em} \eta(\xi)\ :=\ |\xi|J(|\xi|).$$

	It is easy to verify that $\phiep$ satisfies (li1), (li2), (Est),
	(Cpt1), (Cpt2) with
	$$\varphi_{\star}(r)=\varphi^{\star}(r)=r,\hspace{4em}
	\psi_{\star}(r)=\psi^{\star}(r)=r.$$

	Therefore the compactness property (C4) follows from Theorem
	\ref{thm:cpt}.  In order to prove (C1), (C2), (C3) we cannot apply
	directly Theorem \ref{thm:ndim}, since $\varphi$ and $\psi$ do not
	satisfy (\ref{hp:phi}) and (\ref{hp:psi}).  However we can apply
Remark
	\ref{rmk:no-hp-phi}.  Since for $\varphi(r)=\psi(r)=r$ it is well
	known that $\overline{\mathcal{F}_{\varphi,\psi}}(v,\re) =
	|Dv|(\re)$ for every $v\in L^{1}_{loc}(\re)$, by (\ref{ts:rmk}) it
	follows that (C1), (C2), (C3) are satisfied with
	$$\tilde{\mathcal{F}}(u)\ =\
	\int_{\re^n}\left(\int_{\langle\xi\rangle^{\bot}} \left|D
	u_{\xi,y}\right|(\re)\,dy\right)\eta(\xi)\,d\xi\ =\
	\lambda\left|Du\right|(\re^n),$$
	where the last equality follows from a standard integral geometric
	computation.
\end{em}
\end{example}

\begin{example} \label{ex:ms}
\begin{em}
	Let us consider the functionals
	$$\mathcal{F}_{\ep}(u)\ :=\
	\frac{1}{\ep} \int_{\re^n\times\re^n}
	\arctan\left(\frac{\left|u(x+\ep\xi)-u(x)
	\right|^{2}}{\ep|\xi|}\right)J(|\xi|)\,d\xi\,dx.$$

	Then $\left\{\mathcal{F}_{\ep}\right\}$ satisfies (C1) through
	(C4) of \S\ \ref{sec:introduction} with
	$$\mathcal{F}(u)\ :=\
	\cases{{\displaystyle \lambda\int_{\re^n} |\nabla u|^{2} d x +
	\mu\,{\cal H}^{n-1}(S_{u})} & if $u \in \mathit{GSBV}(\re^n)\cap
	L^{1}_{loc}(\re^n)$, \cr \noalign{\vs}
	 +\infty & if $u \in L^{1}_{loc}(\re^n) \setminus
\mathit{GSBV}(\re^n),$ \cr}$$
	where $\lambda=c_{2,n}\,j_{n+1}$, and
	$\mu=\frac{\pi}{2}\,c_{0,n}\,j_{n+1}$.  This family is very
	similar to the family $\left\{\mathcal{DG}_{\ep}\right\}$ studied
	in \cite{gobbino:ms}.  In this case we have that
	$$\phiep(r)\ :=\
	\frac{1}{\ep}\,\arctan(\ep r^{2}),\hspace{4em}
	\eta(\xi)\ :=\ |\xi|J(|\xi|).$$

	Since $\phiep$ satisfies (li1), (li2), (Est), (Cpt1), (Cpt2) with
	$$\varphi_{\star}(r)=\varphi^{\star}(r)=|r|^2,\hspace{4em}
	\psi_{\star}(r)=\psi^{\star}(r)=\frac{\pi}{2},$$
	then the results follow from Theorem \ref{thm:ndim} and Theorem
	\ref{thm:cpt}.
\end{em}
\end{example}

\begin{example} \label{ex:main}
\begin{em}
	Let us consider the functionals
	$$\mathcal{F}_{\ep}(u)\ :=\
	\frac{1}{\ep} \int_{\re^n\times\re^n}
	\frac{\left|u(x+\ep\xi)-u(x)\right|^{2}}{\left
	|u(x+\ep\xi)-u(x)\right|^{3/2}+\ep|\xi|}\,
	J(|\xi|)\,d\xi\,dx.$$

	Then $\left\{\mathcal{F}_{\ep}\right\}$ satisfies (C1) through
	(C4) of \S\ \ref{sec:introduction} with
	$$\mathcal{F}(u)\ :=\
	\cases{{\displaystyle \lambda\int_{\re^n} |\nabla u(x)|^{2}\,dx}+
	{\displaystyle \mu\int_{S_{u}} \sqrt{|u^{+} - u^{-}|}\,d\hn}&
	if $u \in \mathit{GSBV}(\re^n)$, \cr
	\noalign{\vs} +\infty & otherwise, \cr}$$
	where $\lambda=c_{2,n}\ j_{n+1}$, and $\mu=c_{0,n}\,j_{n+1}$.

	Indeed the family $\left\{\mathcal{F}_{\ep}\right\}$ is a
	particular instance of (\ref{defn:fep-ndim-bis}) with
	$$\phiep(r)\
	:=\ \frac{|r|^2}{\sqrt{\ep}|r|^{3/2}+1},\hspace{4em} \eta(\xi)\
	:=\ |\xi|J(|\xi|).$$

	Since it can be proved (exercise for the interested reader!) that
	$\phiep$ satisfies (li1), (li2), (Est) with
	$$\varphi_{\star}(r)=\varphi^{\star}(r)=|r|^2,\hspace{4em}
	\psi_{\star}(r)=\psi^{\star}(r)=\sqrt{r},$$
	then (C1), (C2), (C3) follow from Theorem \ref{thm:ndim}.

	If we want to prove (C4) applying directly Theorem \ref{thm:cpt}, we
	are forced to show that $\phiep$ satisfies (Cpt2), but in this case
	this leads to a huge inequality which seems difficult to prove (or
	disprove). However, we can pursue a different path. We first remark
	that
	\begin{equation}
		\phiep(r)\ \geq\ \tilde{\varphi}_{\ep}(r)\ :=\
		\frac{1}{9}\min\left\{l^{2}+\frac{\sqrt{r-l}}{\sqrt{\ep}} \tc
		l\in[0,r]\right\}
		\label{eqn:mainex}
	\end{equation}
	(use $l=r$ if $r\leq 4\ep^{-1/3}$, and $l=0$ otherwise), and that
	$\{\tilde{\varphi}_{\ep}\}$ satisfies (Cpt1) and (Cpt2) (use Lemma
	\ref{lemma:phiep}). Then we consider the functional
	$\tilde{\mathcal{F}}_{\ep}$ defined as in (\ref{defn:fep-omega})
	with $\tilde{\varphi}_{\ep}$ instead of $\varphi_{\ep}$. Since by
	(\ref{eqn:mainex}) we have that
	$$\sup_{\ep>0}\,\{{\cal
	F}_{\ep}(u_{\ep})+\|u_{\ep}\|_{\infty}\}<+\infty\ \Longrightarrow\
	\sup_{\ep>0}\,\{\tilde{{\cal
	F}}_{\ep}(u_{\ep})+\|u_{\ep}\|_{\infty}\}<+\infty,$$
	then property (C4) for $\left\{\mathcal{F}_{\ep}\right\}$ follows
	applying Theorem \ref{thm:cpt} to the family
	$\{\tilde{\mathcal{F}}_{\ep}\}$.
\end{em}
\end{example}

\end{document}